\documentclass[reqno]{article}
\usepackage{a4wide}
\usepackage{fnpct}
\usepackage[backend=biber,style=alphabetic,backref=true]{biblatex}
\usepackage{graphicx}
\usepackage{xurl}
\usepackage{amsmath,amsthm,amssymb,enumerate,mathtools,environ}
\usepackage[normalem]{ulem} 
\usepackage{euscript,mathrsfs}
\usepackage[left=3cm,right=3cm,top=2.7cm,bottom=2.7cm]{geometry}
\usepackage{color}
\usepackage{upgreek}
\usepackage{bm}
\usepackage{bbm}
\usepackage{parskip}
\usepackage[breaklinks,colorlinks]{hyperref}
\hypersetup{linkcolor=blue,urlcolor=mygreen,citecolor=mygreen}
\catcode`\@=11 \@addtoreset{equation}{section}

\catcode`\@=12
\allowdisplaybreaks
\usepackage[old]{old-arrows}
\newlength{\hookwidth}
\newlength{\hookheight}
\newcommand{\compemb}{%
	\mathrel{%
		\raisebox{-.35\hookheight}{$\varhookrightarrow$}\hspace*{-\hookwidth}%
		\raisebox{ .55\hookheight}{$\varhookrightarrow$}%
	}%
}

\newtheorem{Theorem}{Theorem}[section]
\newtheorem{Proposition}[Theorem]{Proposition}
\newtheorem{Lemma}[Theorem]{Lemma}
\newtheorem{Corollary}[Theorem]{Corollary}
\theoremstyle{definition}
\newtheorem{Definition}[Theorem]{Definition}
\newcommand{\comment}[1]{}

\newtheorem{Remark}[Theorem]{Remark}

\newcommand{\bTheorem}[1]{
\begin{Theorem} \label{T#1} }
\newcommand{\eT}{\end{Theorem}}

\newcommand{\bProposition}[1]{
\begin{Proposition} \label{P#1}}
\newcommand{\eP}{\end{Proposition}}

\newcommand{\bLemma}[1]{
\begin{Lemma} \label{L#1} }
\newcommand{\eL}{\end{Lemma}}

\newcommand{\bCorollary}[1]{
\begin{Corollary} \label{C#1} }
\newcommand{\eC}{\end{Corollary}}

\newcommand{\bRemark}[1]{
\begin{Remark} \label{R#1} }
\newcommand{\eR}{\end{Remark}}

\newcommand{\bDefinition}[1]{
\begin{Definition} \label{D#1} }
\newcommand{\eD}{\end{Definition}}

\newcommand{\dif}{\mathrm{d}}

\newcommand{\mf}{\mathfrak{F}}
\newcommand{\mr}{\mathbb{R}}
\newcommand{\p}{\mathbb{P}}
\newcommand{\T}{\mathbb{T}}
\newcommand{\stred}{\mathbb{E}}

\newcommand{\mn}{\mathbb{N}}

\newcommand{\mt}{\T^N}

\DeclareMathOperator{\supp}{supp}

\newcommand{\bu}{\mathbf u}

\newcommand{\bq}{\mathbf q}

\newcommand{\tor}{\T^N}
\newcommand{\StoB}{\left(\Omega, \mathfrak{F},(\mathfrak{F}_t )_{t \geq 0},  \mathbb{P}\right)}

\newcommand{\bfu}{\mathbf{u}}
\newcommand{\bfv}{\mathbf{v}}
\newcommand{\bfq}{\mathbf{q}}

\newcommand{\bfr}{\mathbf{r}}

\DeclareMathOperator*{\esssup}{ess\,sup}

\newcommand{\bfvarphi}{\boldsymbol{\varphi}}

\newcommand{\bFormula}[1]{
\begin{equation} \label{#1}}
\newcommand{\eF}{\end{equation}}

\newcommand{\vr}{\varrho}

\newcommand{\vu}{\vc{u}}
\newcommand{\vc}[1]{{\bf #1}}

\newcommand{\Div}{{\rm div}_x}
\newcommand{\Grad}{\nabla_x}

\newcommand{\tn}[1]{\mathbb{#1}}
\newcommand{\dx}{\,{\rm d} {x}}

\newcommand{\dt}{\,{\rm d} t }

\newcommand{\vU}{\vc{U}}

\newcommand{\vm}{\vc{m}}

\newcommand{\vv}{\vc{v}}

\newcommand{\D}{{\rm d}}

\newcommand{\R}{\mathbb{R}}

\newcommand{\E}{\mathbb{E}}
\newcommand{\expe}[1]{ \mathbb{E} \left[ #1 \right] }

\definecolor{Cgrey}{rgb}{0.85,0.85,0.85}
\definecolor{Cblue}{rgb}{0.50,0.85,0.85}
\definecolor{Cred}{rgb}{1,.2,.4}
\definecolor{fancy}{rgb}{0.10,0.85,0.10}
\definecolor{mygreen}{rgb}{0.01,0.6,0.2}

\newcommand{\Dif}{{\rm d}}
\newcommand{\DD}{D^d_t}
\newcommand{\DS}{\mathbb{D}^s_t}

\newcommand{\N}{\mathbb N}
\newcommand{\intTorN}[1]{\int_{\T^N} #1 \,{\rm d}x}

\newcommand\Cbox[2]{%
    \newbox\contentbox%
    \newbox\bkgdbox%
    \setbox\contentbox\hbox to \hsize{%
        \vtop{
            \kern\columnsep
            \hbox to \hsize{%
                \kern\columnsep%
                \advance\hsize by -2\columnsep%
                \setlength{\textwidth}{\hsize}%
                \vbox{
                     skip=\baselineskip
                     indent=0bp
                    #2
                }%
                \kern\columnsep%
            }%
            \kern\columnsep%
        }%
    }%
    \setbox\bkgdbox\vbox{
        \color{#1}
        \hrule width  \wd\contentbox %
               height \ht\contentbox %
               depth  \dp\contentbox
        \color{black}
    }%
    \wd\bkgdbox=0bp%
    \vbox{\hbox to \hsize{\box\bkgdbox\box\contentbox}}%
    \vskip\baselineskip%
}

\DeclareMathOperator{\divv}{div}
\DeclareMathOperator{\diver}{div}

\begin{document}
\title{Dissipative Measure Valued Solutions to the Stochastic Compressible Navier--Stokes Equations and Inviscid--Incompressible Limit}

\author{\textsc{Utsab Sarkar}\footnote{{\texttt{utsab@math.iitb.ac.in/reachutsab@gmail.com}}, \textsf{Department of Mathematics, Indian Institute of Technology, Bombay.}}}

\date{}

\maketitle

\begin{abstract}
We introduce the concept of dissipative measure valued martingale solutions for the stochastic compressible Navier--Stokes equations. These solutions are probabilistically weak, as they incorporate both the underlying Wiener process and the probability space as intrinsic components of the formulation. For the stochastic compressible Navier--Stokes system, we further derive a relative energy inequality, which serves as the foundation for establishing a path-wise weak–strong uniqueness principle. We also look at the inviscid--incompressible limit of the underlying system of equations using the relative energy inequality.
\end{abstract}

\providecommand{\keywords}[1]{{\textbf{Keywords :}} #1}

\tableofcontents

\section{Introduction}\label{introduction}
The theory of measure valued solutions has emerged as a powerful framework for addressing the analytical challenges posed by nonlinear partial differential equations (PDEs), particularly in the context of systems exhibiting oscillations or concentrations. The seminal works of DiPerna and Majda \cite{DL,DL1} initiated the study of measure valued solutions for deterministic PDEs associated with conservation laws, introducing a generalized solution concept capable of capturing weak limits of sequences of approximate solutions.

In the setting of compressible fluid dynamics, Neustupa \cite{Neus} was the first to formulate a measure valued solution concept for the compressible Navier--Stokes system. This line of research was further advanced by Feireisl et al. \cite{Fei01}, who introduced the notion of dissipative measure valued solutions, incorporating a relative energy inequality to ensure compatibility with the underlying physical dissipation mechanisms. In this work, we extend the measure valued framework to the stochastic setting by formulating a notion of {\it dissipative measure valued martingale solutions} for the stochastically forced compressible Navier--Stokes system.

We consider a stochastic variant of the {compressible barotropic
Navier--Stokes system} describing the time evolution of the mass density $\vr$ and the bulk velocity $\vu$ of a fluid driven by
a nonlinear multiplicative noise $\mathbb{G}$. The system of equations reads
\begin{align} \label{stochastic ns density equation}
\D \vr + \Div (\vr \vu) \dt &= 0,\\ \label{stochastic ns momentum equation}
\D (\vr \vu) + \left[ \Div (\vr \vu \otimes \vu)+  \Grad p(\vr) \right]  \dt  &= \Div \mathbb{S} (\Grad \vu) \,\dt  + \mathbb{G} (\vr, \vr \vu) \,\D W.
\end{align}
Here
$\mathbb{S}(\Grad \vu)$ is the viscous stress tensor for which we assume Newton's rheological law
\begin{equation} \label{stress tensor}
\mathbb{S}(\Grad \vu) = \nu \left( \Grad \vu + \Grad^t \vu - \frac{2}{N} \Div \vu \mathbb{I} \right) + \lambda \,\Div \vu \mathbb{I}, \qquad \nu > 0, \ \lambda \geq 0.
\end{equation}
Here we consider the pressure to follow the power law, i.e., $p(\varrho)=a\varrho^\gamma$, $\gamma>1$ denotes the adiabatic exponent, $a>0$ is the squared reciprocal of the Mach number (the ratio between average velocity and speed of sound) and to avoid the well-known issues regarding the behaviour of the fluid near { a} boundary and yet considering a physically meaningful situation, we study \eqref{stochastic ns density equation}--\eqref{stochastic ns momentum equation} on the $N$-dimensional flat torus $\T^N$. The driving process $W$ is a cylindrical Wiener process defined on some probability space $(\Omega,\mathfrak{F},\p)$ and the coefficient $\mathbb{G}$ is generally nonlinear and satisfies suitable growth assumptions, see Section \ref{preliminaries and main results} for the precise set-up. The problem is complemented with random initial data
\begin{equation} \label{stochastic nse initial data}
	\vr(0,\cdot) = \vr_0, \ \vr \vu(0, \cdot) = (\vr \vu)_0,
\end{equation}
with sufficient spatial regularity to be specified later. {While the primary focus is on the three-dimensional case $N=3$, our theory also encompasses the cases $N=1\mbox{ and }N=2$, as well as higher (albeit non-physical) dimensions.}

The first existence results were based on an appropriate transformation formula that enables the problem to be reduced to a random system of PDEs where the stochastic integral no longer appears and deterministic methods are available, as shown in \cite{tornatore1998one} for the 1D case and \cite{tornatore2000global} for a very specific periodic 2D case based on the existence theory established by Vaĭgant and Kazhikhov in \cite{vaigant1995existence}. Finally the study by Feireisl et al. \cite{FEIREISL20131342} addresses the 3D scenario. Breit and Hofmanová \cite{BrHo} obtained the first ``really'' stochastic existence result for the compressible Navier--Stokes system perturbed by a general nonlinear multiplicative noise in 3D with periodic boundary conditions where the existence of so-called finite energy weak martingale solutions was proven. Extension of this result to the zero Dirichlet boundary conditions then appeared in \cite{smith2017random,wang2015global}. For completeness, let us also mention \cite{breit2016incompressible}, where a result on singular limit was proved.

Our objective in this work is to formulate an appropriate concept of {\it dissipative measure valued solutions} to the stochastic compressible Navier--Stokes system, as given by equations \eqref{stochastic ns density equation}--\eqref{stochastic ns momentum equation}, and to investigate their fundamental properties. Broadly speaking, these solutions represent a measure valued formulation of the stochastic compressible Navier--Stokes equations that satisfy a suitable form of energy inequality and exist globally in time for arbitrary finite energy initial data. {The motivation for considering measure valued solutions is twofold. First, such solutions provide a natural and potentially maximal class within which the family of smooth (classical) solutions is stable. This is a relatively direct consequence of the weak(measure valued)--strong uniqueness principle. Second, many approximate solutions—such as those arising from numerical schemes or regularized systems—can be shown to converge to a measure valued solution, whereas convergence to weak solutions either remains unknown or involves considerable technical difficulty.} Measure-valued solutions for fluid systems influenced by {\it multiplicative stochastic forcing} constitute a relatively new and developing area of pursuit. To the best of our knowledge, the present work is the first to introduce and analyse a notion of measure valued solution for the stochastic compressible Navier--Stokes equations. We also note that related results concerning dissipative measure valued solutions for the stochastic compressible Euler system were recently established in \cite{mku}.

To motivate our definition of measure valued solution, we consider viscous compressible fluids subject to stochastic forcing described by the following approximate system
\begin{align} \label{P1NS}
\D \vr_{\delta} + \Div (\vr_{\delta} \vu_{\delta}) \,\dt &= 0,\\
 \label{P2NS}
\D (\vr_{\delta} \vu_{\delta}) + \left[ \Div (\vr_{\delta} \vu_{\delta} \otimes \vu_{\delta}) +  \Grad p_{\delta}(\vr) \right] \,\dt  &= \Div \mathbb{S} (\Grad \vu_{\delta}) \,\dt + \mathbb{G} (\vr_{\delta}, \vr_{\delta} \vu_{\delta}) \,\D W.
\end{align}
Here
$p_{\delta}(\vr) = p(\vr) + \delta (\vr + \vr^{\Gamma})$, with $\Gamma \ge \max \lbrace 6, \gamma \rbrace$. In general, the artificial pressure $\delta (\vr + \vr^{\Gamma})$ is required to extract higher integrability of $\vr$ and also to exhibit a better control near vacuum regions.

The study of the system  \eqref{P1NS}--\eqref{P2NS} was first initiated in \cite{BrHo}, where the global-in-time existence of {finite energy weak martingale solutions} is shown. These solutions are weak in both the analytical sense (derivatives only exist in the sense of distributions) as well as the probabilistic sense (the probability space is an integral part of the solution). Moreover, the time-evolution of the energy can be controlled in terms of its initial state.

A central objective of this paper is to establish the weak(measure valued)–strong uniqueness principle for dissipative measure valued martingale solutions to the stochastic compressible Navier--Stokes system. To this end, we first introduce a suitable notion of dissipative measure valued solutions for the stochastic compressible Navier--Stokes equations, specifically designed to facilitate the derivation of the weak–strong uniqueness property within this framework. {A key analytical difficulty arises in the adaptation of techniques from \cite{BrFeHo2015A}, namely, the identification of the martingale terms appearing in the momentum equation \eqref{P2NS} and the associated energy inequality \eqref{energy inequality of approximate nvs} in the limit of a sequence of approximate solutions. In general, such identification proves to be problematic due to the absence of sufficient regularity and uniform estimates on the approximating sequence. In particular, one cannot directly pass to the limit in the martingale terms using either the results of Debussche et al. \cite{Debussche}, or via the strong/weak continuity properties of linear operators such as the Itô integral between Banach spaces. The only available conclusion in this context is that the limiting object retains the martingale property. Nonetheless, a closer examination of the weak–strong uniqueness argument in \cite{BrFeHo2015A} reveals that the full identification of the limiting martingale term is not necessary. Rather, it suffices to characterize the cross-variation between the martingale component of a dissipative solution and the corresponding smooth strong solution.} This crucial insight enables us to establish the weak(measure valued)–strong uniqueness principle for the stochastic compressible Navier--Stokes system \eqref{stochastic ns density equation}--\eqref{stochastic ns momentum equation}. In fact, this is achieved by appropriately specifying the cross-variation with a smooth process, as formalized in item (k) of Definition \ref{dissipative measure valued martingale solution}.

Last but not the least, in order to demonstrate the relative energy inequality's effectiveness, we use it to determine the inviscid--incompressible limit of stochastic compressible flows. It is important to keep in mind that it is fairly challenging to demonstrate directly that, in the low Mach number regime, the solutions of compressible Navier--Stokes equations converge to the solutions of incompressible Euler equations but because of the weak--strong uniqueness property, it is much simpler to demonstrate that the measure valued solutions of the compressible Navier--Stokes equations converge to some genre of strong solutions of the incompressible Euler equations, as it is demonstrated in Section \ref{section incompressible--inviscid limit}. The first significant contribution by Kleinermann and Majda \cite{MajdaK}, in the realm of low Mach limit, was the incompressible limit for the deterministic counterpart of \eqref{stochastic ns density equation}--\eqref{stochastic ns momentum equation}. Many authors have rigorously justified the incompressible limit for the so-called {well-prepared data} (see, for example, \cite{Asano1987, SCHO2, MESC1, Uka}). All of these authors make the assumption that the solutions of compressible flows are smooth. However, regardless of how smooth and/or small the initial data are, strong solutions of the compressible flows produce singularities in a finite amount of time. In the case of compressible inviscid fluids, the hypothesis on the smoothness of solutions is thus highly constrictive and even inappropriate. Considering the aforementioned discussions, to determine the incompressible limit, we adopt the approach put forth by Feireisl et al. \cite{Fei02}, based on the measure valued solution. {We also note that the weak(measure valued)–strong uniqueness principle, along with the concept of dissipative measure valued solutions, has recently been utilized in the analysis of convergence for various numerical schemes. These tools have proven instrumental in establishing rigorous convergence results in complex settings; see, for example, \cite{LM15, FUKRMT17, FLHS21}, and the references therein.} \vspace{2mm}

The rest of the paper is organized as follows :  we describe all the technical frameworks and state the main results in Section~\ref{preliminaries and main results}. Using a stochastic compactness argument, we then present a proof of the existence of dissipative measure valued martingale solutions of the system \eqref{stochastic ns density equation}--\eqref{stochastic nse initial data} in Section~\ref{main proof of existence}. Section~\ref{weakstrong uniqueness} is devoted on deriving the weak(measure valued)--strong uniqueness principle by using a suitable relative energy inequality. Finally, in Section~\ref{section incompressible--inviscid limit}, we explore another application of weak(measure valued)--strong uniqueness property, that justifies the inviscid--incompressible limit of the stochastic compressible Navier--Stokes system \eqref{stochastic ns density equation}--\eqref{stochastic ns momentum equation} rigorously.

\section{Preliminaries and Main Results} \label{preliminaries and main results}

We fix an arbitrary large time horizon, $T>0$. Throughout the paper, we refer to numerous universal constants by the letter $c$. When a constant may change from line to line, we maintain the same notation as long as it does not affect the core idea.

For a generic set $E$, let $\mathcal{M}_b(E)$ denotes the space of bounded Borel measures on $E$, whose norm is given by the total variation norm of measures. It is the dual space of the space of continuous functions vanishing at infinity i.e., $C_0(E)$, equipped with the supremum norm. Moreover, let $\mathcal{P}(E)$ be the space of probability measures on $E$.

\subsection{Analytic Framework}
Throughout the paper, we suppose that, corresponding to the pressure $p = p(\vr)=a\vr^\gamma$, the pressure potential is, $P(\varrho)=\varrho\int_1^\varrho \frac{p(z)}{z^2}\,\mathrm{d}z$.

For a separable Hilbert space $H$ and $\alpha\in(0,1),{\,p>1}$; let $W^{\alpha,p}(0,T;H)$ be the $H$-valued {fractional} Sobolev space, characterised by the norm \[\|u\|^p_{W^{\alpha,p}(0,T;H)}\coloneqq\int_0^T\|u\|^p_{H}\dt+\int_0^T\int_0^T\frac{\|u(t)-u(s)\|^p_{H}}{|t-s|^{1+p\alpha}}\,\dt\,\D s.\]
We recall here the following compact embedding result form Flandoli and Gatarek \cite[Theorem 2.2]{FlandoliGatarek}.
\begin{Lemma}\label{FlandoliGatarek}
If $X,Y$ are two Banach spaces with the compact embedding $X\compemb Y$, and the real numbers $\alpha\in\left(0,1\right),\,p>1$ satisfy $\alpha p>1$, then the following compact embedding holds. \[W^{\alpha,p}(0,T;X)\compemb C([0,T];Y).\]
\end{Lemma}

\subsubsection{Young Measures and Concentration Defect Measures}
We will briefly review here the concept of Young measures and associated results, which have been referenced extensively in the paper. We cite Balder \cite{Balder} for an excellent summary of Young measure theory applications to hyperbolic conservation laws. In what follows, lets assume first that $(X,\mathcal{M},\mu)$ is a $\sigma$-finite measure space and a Young measure from $X$ into $\mathbb{R}^d$ is a weakly-$*$ measurable map $\nu:X\rightarrow\mathcal{P}(\mathbb{R}^d)$, in the sense that $x\mapsto\nu_x(B)$ is $\mathcal{M}$-measurable for every Borel set $B$ in $\mathbb{R}^d$.

We make use of the following generalization of the classical result on Young measures, see \cite[Section 2.8]{BrFeHobook}.
\begin{Lemma}\label{younglemma}
Let $N,M\in\mathbb{N},\,\mathcal{O}\subset\mathbb{R}^N\times(0,T)$ and let $({\bf U}_n)_{n \in \N}, {\bf U}_n:\Omega\times\mathcal{O}\rightarrow\R^M$, be a sequence of random variables such that
$$
\E \big[ \|{\bf U}_n \|^p_{L^p(\mathcal{O})}\big] <+ \infty, \,\, \text{for a certain}\,\, p\in(1,\infty).
$$
Then there exists a new subsequence $(\tilde {\bf U}_n)_{n \in \N}$ (not relabelled) and a parameterized family ${\lbrace \tilde \nu^\omega_{x,t} \rbrace}_{(x,t) \in \mathcal{O}}$ of random probability measures on $\R^M$, defined on the standard probability space $\big([0,1], \overline{\mathcal{B}[0,1]}, \mathcal{L} \big)$, such that
$$
{\bf U}_n \sim_{d} \tilde {\bf U}_n,
$$
with the property that, for any Carath\'eodory function $G=G(\theta, Z), \theta \in \mathcal{O}, Z \in \R^M$, such that
$$
|G(\theta,Z)| \le c(1 + |Z|^q), \quad 1 \le q < p, \,\,\text{uniformly in}\,\,\theta,
$$
implies $\mathcal{L}$-{\normalfont a.s.},
$$
G(\cdot, \tilde {\bf U}_n) \rightharpoonup \overline{G}\,\,\text{in}\,\, L^{p/q}(\mathcal{O}), 
\,\, \text{where}\,\, \overline{G}(\theta)= \int_{\R^M} G(\theta, z)\,\D \nu_{(t,x,\omega)}(z), \,\,\text{for {\normalfont a.e.}}\,\, \theta \in \mathcal{O}.
$$
\end{Lemma}
 We know that Young measure theory works well when we try to extract limits of bounded continuous functions. Next, we examine what happens for those functions $H$ for which we only know that 
\begin{equation*}
\E\| H(\vU_n)\|^p_{L^1(\mathcal{O})} \leq c, \quad \mbox{uniformly in } n,\text{ for a certain }p\in(1,\infty).
\end{equation*} By using the fact that $L^1(\mathcal{O})$ is embedded in the space of bounded Radon measures $\mathcal{M}_b(\mathcal{O})$, we can infer that $\p$-a.s.
\begin{align*}
	\mbox{weak$^*$ limit in} \, \mathcal{M}_b(\mathcal{O})\,\, \mbox{of} \, \,H({\bf U}_n) = \langle \mathcal{\tilde V}^{\omega}_{x}; H\rangle \,\dx+ H_{\infty},
\end{align*}
where $H_{\infty} \in \mathcal{M}_b(\mathcal{O})$, and $H_{\infty}$ is called {concentration defect measure or concentration Young measure}. Note that, a simple truncation analysis and Fatou's lemma reveal that $\tilde\p$-a.s. $\| \langle \mathcal{\tilde V}^{\omega}_{(\cdot)}; H\rangle\|_{L^1(\mathcal{O})} \leq c$ and thus $\tilde\p$-a.s. $\langle \mathcal{\tilde V}^{\omega}_{x}; H \rangle $ is finite for a.e. $x\in \mathcal{O}$. In what follows, regarding the concentration defect measure, we shall make use of the following crucial lemma. For a proof of this lemma modulo cosmetic changes, we refer to Feireisl et. al. \cite[Lemma 2.1]{Fei01}. 

\begin{Lemma}
	\label{lemma001}
	Let $\{{\bf U}_n\}_{n > 0}$, ${\bf U}_n: \Omega \times \mathcal{O} \rightarrow \mathbb{R}^M$ be a sequence generating a Young measure $\{\mathcal{V}^{\omega}_y\}_{y\in \mathcal{O}}$, where $\mathcal{O}$ is a measurable set in $\mathbb{R}^N \times (0,T)$. Let $G: \mathbb{R}^M \rightarrow [0,\infty)$
	be a continuous function such that
	\begin{equation*}
		\sup_{n \geq 1} \E\big[\|G({\bf U}_n)\|^p_{L^1(\mathcal{O})} \big]< \infty, \, \text{for a certain}\,\, p\in(1,\infty),
	\end{equation*}
	and let $F$ be continuous such that
	\begin{equation*}
		F: \mathbb{R}^M \rightarrow \mathbb{R}, \quad |F(\bm{z})|\leq G(\bm{z}), \mbox{ for all } \bm{z}\in \mathbb{R}^M.
	\end{equation*}
	Let us denote $\p$-{\normalfont a.s.}
	\begin{equation*}
		{F_{\infty}} \coloneqq {\tilde{F}}- \langle \mathcal{ V}^{\omega}_{y}, F \rangle \,\D y, \quad 
		{G_{\infty}} \coloneqq {\tilde{G}}- \langle \mathcal{ V}^{\omega}_{y}, G \rangle \,\D y.
	\end{equation*}
	Here ${\tilde{F}}, {\tilde{G}} \in \mathcal{M}_b(\mathcal{O})$ are weak$^*$ limits of $\{F({\bf U}_n)\}_{n \geq1}$, $\{G({\bf U}_n)\}_{n \geq 1}$ respectively in $\mathcal{M}_b(\mathcal{O})$. Then $\p${\normalfont-a.s.} $|F_{\infty}| \leq G_{\infty}$.
\end{Lemma}

\subsection{Stochastic Framework}

Here we specify the stochastic forcing term.
Let $(\Omega,\mf,(\mf_t)_{t\geq0},\mathbb{P})$ be a stochastic basis with a complete, right-continuous filtration. The process $W$ is a cylindrical Wiener process in a separable Hilbert space $\mathfrak U$ and has the following representation
$$W(t)=\sum_{k\geq 1} e_k W_k(t),$$
where $\{ W_k \}_{k \geq 1}$ is a family of mutually independent real valued Brownian motions and $\{e_k\}_{k\geq 1}$ is an orthonormal basis of $\mathfrak{U}$. To give the precise definition of the diffusion coefficient $\mathbb{G}$, consider $\varrho\in L^\gamma(\mt)$, $\varrho\geq0$, and $\bfv\in L^2(\mt)$ such that $\sqrt\varrho\bfv\in L^2(\mt)$. Let $\bfq=\varrho\bfv$ and $\,\mathbb{G}(\varrho,\bq):\mathfrak{U}\rightarrow L^1(\mt)$ be defined as follows
$$\mathbb{G}(\varrho,\bq)e_k=\mathbf{G}_k(\cdot,\varrho(\cdot),\bq(\cdot)).$$
The coefficients $\mathbf{G}_{k}:\mt\times\mr\times\R^N\rightarrow\R^N$ are $C^1$-functions that satisfy uniformly in $x\in\mt$
\begin{align}
\vc{G}_k (\cdot, 0 , 0) &= 0, \label{noise condition 1}\\
|  \partial_\vr \vc{G}_k | + |\nabla_{\vc{q}} \vc{G}_k | &\leq \alpha_k, \quad \sum_{k \geq 1} \alpha_k  < \infty.
\label{noise condition 2}
\end{align}
As in \cite{BrHo}, we understand the stochastic integral $\mathbb G\D W$ as a process in the Hilbert space $W^{-m,2}(\mt)$ for $m>N/2$. Indeed, it can be checked that under the above assumptions on $\varrho$ and $\bfv$, the mapping $\mathbb{G}(\varrho,\varrho\bfv)$ belongs to $L_2(\mathfrak{U};W^{-m,2}(\mt))$, the space of Hilbert-Schmidt operators from $\mathfrak{U}$ to $W^{-m,2}(\mt)$.
Consequently, if\footnote{Here $\mathcal{P}$ denotes the predictable $\sigma$-algebra associated to $(\mf_t)$.}
\begin{align*}
\varrho&\in L^\gamma(\Omega\times(0,T),\mathcal{P},\dif\mathbb{P}\otimes\dif t;L^\gamma(\mt)),\\
\sqrt\varrho\bfv&\in L^2(\Omega\times(0,T),\mathcal{P},\dif\mathbb{P}\otimes\dif t;L^2(\mt)),
\end{align*}
and the mean value $(\varrho(t))_{\mt}$ is essentially bounded then the stochastic integral
\[
\int_0^t \mathbb{G}(\vr, \vr \vu) \ {\rm d} W = \sum_{k \geq 1}\int_0^t \vc{G}_k (\cdot, \vr, \vr \vu) \ {\rm d} W_k
\]
 is a well-defined $(\mf_t)$-martingale taking values in $W^{-m,2}(\mt)$. Note that the continuity equation \eqref{stochastic ns density equation} implies that the mean value $(\varrho(t))_{\mt}$ of the density $\varrho$ is constant in time (but in general depends on $\omega$).
Finally, we define the auxiliary space $\mathfrak{U}_0\supset\mathfrak{U}$ via
$$\mathfrak{U}_0=\bigg\{v=\sum_{k\geq1}\alpha_k e_k;\;\sum_{k\geq1}\frac{\alpha_k^2}{k^2}<\infty\bigg\},$$
endowed with the norm
$$\|v\|^2_{\mathfrak{U}_0}=\sum_{k\geq1}\frac{\alpha_k^2}{k^2},\quad v=\sum_{k\geq1}\alpha_k e_k.$$
Note that the embedding $\mathfrak{U}\hookrightarrow\mathfrak{U}_0$ is Hilbert-Schmidt. Moreover, trajectories of $W$ are $\mathbb{P}$-a.s. in $C([0,T];\mathfrak{U}_0)$.

Next, we recall a modified version of the classical Skorokhod embedding theorem \cite[Theorem 2]{Jakubowski}, as stated below.
\begin{Theorem}[Jakubowski--Skorokhod]
	\label{JakubowskiSkorohod}
	Let $(\Omega, \mathcal{F},\mathbb{P})$ be a probability space, and $H$ be a quasi-polish space (i.e., there is a sequence of continuous functions $h_n:H\to[-1,1]$ that separates points of $H$). Assume that $\mathcal{H}$ is the sigma algebra generated by the sequence of continuous functions $\{h_n\}_{n=1}^{\infty}$. Let $U_n:\Omega\to H$, $n\in\mathbb{N}$, be a family of random variables, such that $\{\mathcal{L}aw(U_n):n\in\mathbb{N}\}$ is a tight sequence of probability measures on $(H, \mathcal{H})$. Then, there exists a probability space $(\tilde{\Omega},\tilde{\mathcal{F}},\tilde{\mathbb{P}})$, a family of $H$-valued random variables $\{\tilde{U}_n:n\in\mathbb{N}\}$, on $(\tilde{\Omega},\tilde{\mathcal{F}},\tilde{\mathbb{P}})$ and a random variable $\tilde{U}:\tilde{\Omega}\to H$ such that
	\begin{enumerate}
		\item [\normalfont{(a)}] $\mathcal{L}aw(\tilde{U}_n)=\mathcal{L}aw(U_n),\,\forall\,n\in\mathbb{N};$
		\item [\normalfont{(b)}]
		$\tilde{U}_n\to \tilde{U}\,\text{in}\,H,\,\tilde\p$-{\normalfont a.s.} Moreover, the law of $\tilde{U}$ is a Radon measure.
	\end{enumerate}
\end{Theorem}

Finally, we recall the ``Kolmogorov test'' for the existence of continuous modifications of real valued stochastic processes.

\begin{Lemma}\label{lemma01}
	Let $Y={\lbrace Y(t) \rbrace}_{t \in [0,T]} $ be a real valued stochastic process defined on a probability space $(\Omega,\mf,(\mf_t)_{t\geq0},\mathbb{P})$. Suppose that there are constants $\alpha >1, \beta >0$, and $C>0$ such that for all $s,t \in [0,T]$,
	\begin{align*}
		\E[|Y(t)-Y(s)|^\alpha] \le C |t-s|^{1 + \beta}.
	\end{align*}
	Then there exists a continuous modification of $Y$ and the paths of $Y$ are $\kappa$-H\"{o}lder continuous for every $\kappa \in[0, \frac{\beta}{\alpha})$.
\end{Lemma}

\subsection{Stochastic Compressible Navier--Stokes System}
We introduce the notions of local {and maximal} strong path-wise solutions. Such solutions are strong in both the analytical and the probabilistic sense but exists only locally in time. To be more precise, system \eqref{stochastic ns density equation}--\eqref{stochastic ns momentum equation} will be satisfied pointwise (not only in the sense of distributions) on the given stochastic basis associated to the cylindrical Wiener process $W$.

\begin{Definition}[Local Strong Path-wise Solution] \label{def strong pathwise solution for nse}
Let $\StoB$ be a stochastic basis with a complete right-continuous filtration. Let ${W}$ be an $(\mathfrak{F}_t) $-cylindrical Wiener process and $(\varrho_0,\bfu_0)$ be a $W^{s,2}(\T^N)\times W^{s,2}(\T^N;\R^N)$-valued $\mathfrak{F}_0$-measurable random variable for some ${s>\frac{N}{2}+2}$. A triplet
$(\varrho,\vu,\mathfrak{t})$ is called a local strong path-wise solution to the system \eqref{stochastic ns density equation}--\eqref{stochastic nse initial data} provided
\begin{enumerate}[(a)]
\item $\mathfrak{t}$ is an a.s. strictly positive  $(\mathfrak{F}_t)$-stopping time;
\item the density $\varrho$ is a $W^{s,2}(\mt)$-valued $(\mathfrak{F}_t)$-progressively measurable process satisfying
$$\varrho(\cdot\wedge \mathfrak{t})  > 0,\ \varrho(\cdot\wedge \mathfrak{t}) \in C([0,T]; W^{s,2}(\mt)) \quad \mathbb{P}\text{-a.s.};$$
\item the velocity $\vu$ is a $W^{s,2}(\mt;\R^N)$-valued $(\mathfrak{F}_t)$-progressively measurable process satisfying
$$ \vu(\cdot\wedge \mathfrak{t}) \in   C([0,T]; W^{s,2}(\mt)) \cap L^2(0,T; W^{s+1,2}(\T^N))\quad \mathbb{P}\text{-a.s.};$$
\item  there holds $\mathbb{P}$-a.s.
\[
\begin{split}
&\varrho (t\wedge \mathfrak{t}) = \varrho_0 -  \int_0^{t \wedge \mathfrak{t}} \Div(\varrho\vu ) \ \dif s, \\
&(\varrho \vu) (t \wedge \mathfrak{t})  = \varrho_0 \vu_0 - \int_0^{t \wedge \mathfrak{t}} \Div (\varrho\vu \otimes\vu ) \ \dif s \\
& \hspace{2cm} - \int_0^{t \wedge \mathfrak{t}} \Grad p(\varrho) \ \dif s +\int_0^{t \wedge \mathfrak{t}} \Div \mathbb{S} (\Grad \vu) \ \dif s + \int_0^{t \wedge \mathfrak{t}} {\tn{G}}(\varrho,\varrho\vu ) \ \D W,
\end{split}
\]
for all $t\in[0,T]$.
\end{enumerate}
\end{Definition}

In the above definition, we have assumed $s$ is large enough in order to provide sufficient regularity for the strong solutions. Classical solutions require {two} spatial derivatives of $\vu$ to be continuous $\mathbb{P}$-a.s. This motivates the following definition.

\begin{Definition}[Maximal Strong Path-wise Solution]\label{def maximal strong pathwise solution for nse}
Fix a stochastic basis with a cylindrical Wiener process and an initial condition as in Definition \ref{def strong pathwise solution for nse}. A quadruplet $$(\varrho,\vu,(\mathfrak{t}_R)_{R\in\mn},\mathfrak{t})$$ is a maximal strong path-wise solution to system \eqref{stochastic ns density equation}--\eqref{stochastic nse initial data} provided

\begin{enumerate}[(a)]
\item $\mathfrak{t}$ is an a.s. strictly positive $(\mathfrak{F}_t)$-stopping time;
\item $(\mathfrak{t}_R)_{R\in\mn}$ is an increasing sequence of $(\mathfrak{F}_t)$-stopping times such that
$\mathfrak{t}_R<\mathfrak{t}$ on the set $[\mathfrak{t}<T]$, $\lim_{R\to\infty}\mathfrak{t}_R=\mathfrak t$ a.s. and
\begin{equation}\label{eq:blowup}
\sup_{t\in[0,\mathfrak{t}_R]}\|\vu(t)\|_{W^{2,\infty}}\geq R\quad \text{on}\quad [\mathfrak{t}<T] ;
\end{equation}
\item each triplet $(\varrho,\vu,\mathfrak{t}_R)$, $R\in\mn$,  is a local strong path-wise solution in the sense  of Definition \ref{def strong pathwise solution for nse}.
\end{enumerate}
\end{Definition}

The notion of a maximal path-wise solution has already appeared in the literature in the context of various SPDE or SDE models, see for instance \cite{BMS,Elw, jac,MikRoz}. We state here the main result from \cite[Theorem 5.0.3]{BrFeHobook} or see \cite[Theorem 2.7]{BrFeHo2016}.

\begin{Theorem}\label{existence of finite enerey weak martingale solution for nse}
Let $s\in\mn$ satisfy $s>\frac{N}{2}+3$. Let the coefficients
$\mathbf{G}_k$ satisfy hypotheses \eqref{noise condition 1}, \eqref{noise condition 2} with $C^s$-regularity and let $(\varrho_0,\bfu_0)$ be an $\mathfrak{F}_0$-measurable, $W^{s,2}(\mt)\times W^{s,2}(\mt;\R^N)$-valued random variable such that $\varrho_0>0$ $\p${\normalfont-a.s}. Then there exists a unique maximal strong path-wise solution $(\varrho,\vu,(\mathfrak{t}_R)_{R\in\mn},\mathfrak{t})$ to problem \eqref{stochastic ns density equation}--\eqref{stochastic nse initial data}
in the sense of Definition {\normalfont\ref{def maximal strong pathwise solution for nse}} with the initial condition $(\varrho_0,\vu_0)$.
\end{Theorem}

\subsection{Approximate Compressible Navier--Stokes System}
In this section, we recall the concept of finite energy weak martingale solutions
to \eqref{P1NS}--\eqref{P2NS}, that has been introduced in \cite{BrHo}
and improved in \cite{BrFeHo2015A} and \cite{BrFeHobook}. Both papers complement \eqref{P1NS}--\eqref{P2NS} with periodic boundary conditions. A corresponding version on the whole space can be found in \cite{Romeo1}. These solutions are weak in both the analytical sense and the probabilistic sense. Moreover, the time-evolution of the energy can be controlled in terms of its initial state and these solutions exist globally in time.

\begin{Definition}[Finite Energy Weak Martingale Solution]
\label{finite energy weak martingale solution}
Let $\Lambda_\delta$ be a Borel probability measure on $L^1(\mathbb{T}^N)\times L^1(\mathbb{T}^N)$. Then $\big[ \big(\Omega_{\delta},\mathfrak{F}_{\delta}, (\mathfrak{F}_{{\delta},t})_{t\geq0},\mathbb{P}_{\delta} \big);\varrho_{\delta},\mathbf{u}_{\delta}, W_{\delta} \big]$ is called a finite energy weak martingale solution to \eqref{P1NS}--\eqref{P2NS} with the initial law $\Lambda_\delta$, provided :

\begin{enumerate}[(a)]
	\item $\big(\Omega_{\delta},\mathfrak{F}_{\delta}, (\mathfrak{F}_{{\delta},t})_{t\geq0},\mathbb{P}_{\delta} \big)$ is a stochastic basis with a complete right-continuous filtration;
	\item $W_\delta$ is a cylindrical $(\mathfrak{F}_{\delta,t})$-Wiener process;
	\item the density $\varrho_{\delta}$ and the velocity $\mathbf{u}_\delta$ are random distributions adapted to $(\mathfrak{F}_{\delta,t})$, $\varrho_{\delta}\geq 0$ $\mathbb{P}_\delta$-a.s.;
	\item there exists an $\mathfrak{F}_0$-measurable random variable $[\varrho_{\delta,0},\mathbf{u}_{\delta,0}]$ such that \[\Lambda_\delta=\mathcal{L}[\varrho_{\delta,0},\varrho_{\delta,0}\mathbf{u}_{\delta,0}].\]
	\item the equation of continuity \[-\int_0^T \partial_t\phi\int_{\mathbb{T}^N}\varrho_{\delta}\psi\ \mathrm dx\ \mathrm dt=\phi(0)\int_{\T^N}\varrho_{\delta,0}\psi\ \mathrm{d}x+\int_0^T\phi\int_{\T^N}\varrho_{\delta}\mathbf{u}_\delta\cdot\nabla_x\psi\ \mathrm{d}x\ \mathrm dt\] holds for all $\phi\in C_c^\infty[0,T)$ and all $\psi\in C^\infty(\T^N)$ $\mathbb{P}_\delta$-a.s.;
	\item the equation of momentum
	\begin{align*}
	&-\int_0^T \partial_t\phi\int_{\mathbb{T}^N}\varrho_{\delta}\mathbf{u}_\delta\cdot\upvarphi\ \mathrm dx\ \mathrm dt-\phi(0)\int_{\T^N}\varrho_{\delta,0}\mathbf{u}_{\delta,0}\cdot\upvarphi\ \mathrm dx\\ &\ \ =\int_0^T\phi\int_{\T^N}\left[\varrho_{\delta}\mathbf{u}_\delta\otimes\mathbf{u}_\delta\colon\nabla_x\upvarphi+p_\delta(\varrho_{\delta})\divv_x\upvarphi\right]\mathrm dx\ \mathrm dt\\ &\qquad -\int_0^T\phi\int_{\T^N}\mathbb{S}(\nabla_x\mathbf{u}_\delta)\colon\nabla_x\upvarphi\ \mathrm dx\ \mathrm dt+\int_0^T\phi\int_{\T^N}\mathbb{G}(\varrho_{\delta},\varrho_{\delta}\mathbf{u}_\delta)\cdot\upvarphi\ \mathrm dx\ \mathrm dW_\delta 
	\end{align*} holds for all $\phi\in C_c^\infty[0,T)$ and all $\upvarphi\in C^\infty(\T^N;\R^N)$ $\mathbb{P}_\delta$-a.s.;
	\item the energy inequality \begin{equation}\label{energy inequality of approximate nvs}
		\begin{aligned}
	&-\int_0^T \partial_t\phi\int_{\T^N}\left[\frac12\varrho_{\delta}|\mathbf u_\delta|^2+P_\delta(\varrho_{\delta})\right]\mathrm dx\ \mathrm dt-\phi(0)\int_{\T^N}\left[\frac12\varrho_{\delta,0}|\mathbf u_{\delta,0}|^2+P_\delta(\varrho_{\delta,0})\right]\mathrm dx\\
	&\ \ + \int_0^T\phi\int_{\T^N}\mathbb{S}(\nabla_x\mathbf{u}_\delta)\colon\nabla_x\mathbf{u}_\delta\ \mathrm dx\ \mathrm dt\\
	&\leq \frac12 \sum_{k=1}^\infty\int_0^T\phi\int_{\T^N}\varrho_{\delta}^{-1}|\mathbf G_k(\varrho_{\delta},\varrho_{\delta}\mathbf{u}_\delta)|^2\ \mathrm dx\ \mathrm dt+\int_0^T\phi\int_{\T^N}\mathbb{G}(\varrho_{\delta},\varrho_{\delta}\mathbf{u}_\delta)\cdot\mathbf{u}_\delta\ \mathrm dx\ \mathrm dW_\delta
\end{aligned}
	\end{equation}
	holds for all $\phi\in C_c^\infty[0,T)$, $\phi\geq0$, $\mathbb{P}_\delta$-a.s.
\end{enumerate}
\end{Definition}	

\begin{Remark}
	Note that we may take, for every $\delta$, $\Omega_{\delta} =[0,1]$, with $\mathfrak{F}_{\delta}$ the $\sigma$-algebra of the Borelians on $[0,1]$, and $\mathbb{P}_{\delta}$ the Lebesgue measure on $[0,1]$, see Skorokhod \cite{sko}. Moreover, we can assume without loss of generality that there exists one common Wiener space $W$ for all $\delta$. Indeed, this can be achieved by performing the compactness argument from any chosen subsequence ${\lbrace \delta_n \rbrace}_{n \in \N}$ at once. However, we stress on the fact that it may not be possible to exclude the dependence of $\delta$ in the filtration $ (\mathfrak{F}_{{\delta},t})_{t\geq0}$, due to lack of path-wise uniqueness for the underlying system.
\end{Remark}\vspace{2mm}

The following existence theorem is shown in \cite[Theorem 4.4.2]{BrFeHobook}.
\begin{Theorem}
\label{thm:BrFeHobook}
Let $\Gamma\geq6$ and assume that $\Lambda_{\delta}$ is a Borel probability measure on $L^1(\mathbb{T}^3)\times L^{1}(\mathbb{T}^3)$ such that
\begin{align*}
\Lambda_{\delta}\Bigg\{(\varrho,\mathbf{q})\in L^1(\mathbb{T}^3)\times L^{1}(\mathbb{T}^3) \, :\, 0<M_1\leq \int_{\T^3}\varrho\ \mathrm dx\leq M_2,\, \varrho>0 \Bigg\}=1, 
\end{align*}
with deterministic constants $0<M_1<M_2<\infty$. Furthermore, assume that the following moment estimate
\begin{align}\label{initial}
\int_{L^1_{x}\times  L^{1}_x} \bigg\Vert\frac{1}{2}\frac{\vert  \mathbf{q} \vert^2}{\varrho}+ P_{\delta}(\varrho)  \bigg\Vert^r_{L^1_x(\T^3)}\, \mathrm{d}\Lambda_{\delta}(\varrho,\mathbf{q})<\infty
\end{align}
holds for some $r\geq4$. Finally, also assume \eqref{noise condition 1}, \eqref{noise condition 2}. Then there exists a finite energy weak martingale solution of \eqref{P1NS}--\eqref{P2NS} in the sense of Definition {\normalfont\ref{finite energy weak martingale solution}} with initial law $\Lambda_{\delta}$.
\end{Theorem}

{We remark that the most intricate and physically significant scenario arises when \( N = 3 \). The proof of Theorem~\ref{thm:BrFeHobook} extends to any dimension \( N \ge 1 \) with straightforward modifications. For further details, we direct interested readers to \cite[Chapter~4]{BrFeHobook}.}

\subsection{Measure Valued Solutions}
In general, the energy inequality \eqref{energy inequality of approximate nvs} seems to be the only source of {\it a priori} bounds. However, as indicated by the numerous
examples of ``oscillatory'' solutions the set of all admissible weak solutions emanating from given initial data is not closed with respect to the weak topology on the trajectory space associated with the energy bounds. There are two potential sources of difficulties:

\begin{itemize}

\item non-controllable {oscillations} due to accumulation of singularities;

\item blow-up type collapse due to possible {concentration} points.

\end{itemize}

To accommodate these singularities in the closure of the set of weak solutions, two kinds of tools are used:
(i) the Young measures describing the oscillations, (ii) concentration defect measures for concentrations.

Denote
\[
\mathcal{S} = \left\{ [\vr, \vr\vu] \ \Big| \,\vr \geq 0,\,\vr\vu \in \R^N \right\}
\]
to be the phase space associated to the Navier--Stokes system \eqref{stochastic ns density equation}--\eqref{stochastic nse initial data}. Let $\mathcal{P}(\mathcal{S})$ denote the set of probability measures on $\mathcal{S}$ and let $\mathcal{M}^{+}(\Theta)$ denote the set of positive bounded Radon measures on $\Theta$.


\subsubsection{Dissipative Measure Valued Solutions for the Stochastic Compressible Navier--Stokes System}
	
Motivated by the previous discussion, we are ready to introduce the concept of {dissipative measure valued solution} to the stochastic compressible Navier--Stokes system.

\begin{Definition}[Dissipative Measure Valued Martingale Solution]
\label{dissipative measure valued martingale solution}
Let $\Lambda$ be a Borel probability measure on $L^\gamma(\T^N)\times L^{2\gamma/\gamma+1}(\T^N)$. Then $\big[ \big(\Omega,\mathfrak{F}, (\mathfrak{F}_{t})_{t\geq0},\mathbb{P} \big); \nu_{(t,x,\omega)}, W \big]$ is a dissipative measure valued martingale solution of \eqref{stochastic ns density equation}--\eqref{stochastic ns momentum equation}, with initial condition $\nu_{(0,x,\omega)}$, if
\begin{enumerate}[(a)]
\item $\nu_{(t,x,\omega)}$ is a parametrized family of probability measures
\begin{equation*}
\nu_{(t,x,\omega)}: (t,x,\omega) \in [0,T] \times \T^N \times \Omega \rightarrow \mathcal{P}(\mathcal{S}).
\end{equation*}
That is, $\nu$ can be regarded as a random variable taking values in the space of Young measures, namely $L^{\infty}_{w^*}\big([0,T] \times \T^N; \mathcal{P}\big(\mathcal{S})\big)$.
\item $ \big(\Omega,\mathfrak{F}, (\mathfrak{F}_{t})_{t\geq0},\mathbb{P} \big)$ is a stochastic basis with a complete right-continuous filtration,
\item $W$ is a $(\mathfrak{F}_{t})$-cylindrical Wiener process,
\item the average density $\langle \nu_{(t,x,\omega)}; \varrho \rangle$ satisfies $t\mapsto \langle \langle \nu_{(t,x,\omega)}; \varrho \rangle(t, \cdot),\phi\rangle\in C[0,T]$ for any $\phi\in C^\infty(\T^N)$ $\mathbb{P}$-a.s., the function $t\mapsto \langle \langle \nu_{(t,x,\omega)}; \varrho \rangle(t, \cdot),\phi\rangle$ is progressively measurable and
\begin{align*}
\mathbb{E}\, \bigg[ \sup_{t\in(0,T)}\Vert  \langle \nu_{(t,x,\omega)}; \varrho \rangle(t,\cdot)\Vert_{L^\gamma(\mt)}^p\bigg]<\infty,
\end{align*}
for all $1\leq p<\infty$,
\item the average velocity $\big<\nu_{(t,x,\omega)}; {\bf u} \big>  \in L^2(0,T; W^{1,2}(\T^N;\R^N))$, and the average momentum $\langle \nu_{(t,x,\omega)};\textbf{m} \rangle\newline=\langle \nu_{(t,x,\omega)}; \varrho\textbf{u} \rangle$ satisfies  $t\mapsto \langle \langle \nu_{(t,x,\omega)}; \varrho\textbf{u} \rangle (t, \cdot),\bm\phi\rangle\in C[0,T]$ for any $\bm{\phi}\in C^\infty(\T^N;\R^N)$ $\mathbb{P}$-a.s., the function $t\mapsto \langle \langle \nu_{(t,x,\omega)}; {\bf m} \rangle (t, \cdot),\bm{\phi}\rangle$ is progressively measurable and 
\begin{align*}
\mathbb{E}\, \bigg[ \sup_{t\in(0,T)}\Vert  \langle \nu_{(t,x,\omega)}; \textbf{m} \rangle (t, \cdot) \Vert_{L^\frac{2\gamma}{\gamma+1}(\mt)}^p\bigg]<\infty,
\end{align*}
for all $1\leq p<\infty$,
\item $\Lambda=\mathcal{L}[\nu_{(0,x,\omega)}]$, where $\nu_{(0,x,\omega)}=\delta_{[{\varrho}(0,x,\omega),({\varrho}{\bf u})(0,x,\omega)]}$,
\item the integral identity
\begin{equation} 
\begin{aligned}\label{measure valued density equation}
&\int_{\T^N} \langle \nu_{(\tau,x,\omega)}; \varrho \rangle \, \varphi(x)\,\dx -\int_{\T^N} \langle \nu_{(0,x,\omega)}; \varrho \rangle\, \varphi (x)\,\dx \\
& \qquad \qquad = \int_{0}^{\tau} \int_{\T^N}  \langle \nu_{(t,x,\omega)}; \textbf{m}\rangle \cdot \nabla_x \varphi \,\dx \,\dt + \int_{0}^{\tau}\int_{\T^N} \nabla_x \varphi \cdot \D\mu_c
\end{aligned}
\end{equation}
holds for all $\tau \in[0,T)$, $\p$-a.s., and for all $\varphi \in C^{\infty}(\T^N)$, where $\mu_c \in \mathcal{M}([0,T]\times \T^N; \mathbb{R}^N)$, $\p$-a.s., is a vector valued measure;
\item there exists a $W^{-m,2}(\T^N)$-valued square integrable continuous martingale $M^1_{E}$, such that the integral identity
\begin{equation} \label{measure valued momentum equation}
\begin{aligned}
&\int_{\T^N} \langle \nu_{(\tau,x,\omega)}; \textbf{m}\rangle \cdot \bm{\varphi}(x) \dx - \int_{\T^N} \langle \nu_{(0,x,\omega)}; \textbf{m}\rangle \cdot \bm{\varphi}(x) \dx \\
&\quad = \int_{0}^{\tau} \int_{\T^N} \left[\left\langle \nu_{(t,x,\omega)}; \frac{\textbf{m}\otimes \textbf{m}}{\varrho} \right\rangle: \nabla_x \bm{\varphi} + \langle \nu_{(t,x,\omega)};p(\varrho)\rangle \divv_x \bm{\varphi} \right] \dx\,\dt \\&\quad\quad-  \int_{0}^{\tau} \int_{\T^N} \mathbb{S} (\nabla_x  {\langle\nu_{(t,x,\omega)},\bf u\rangle}) : \nabla_x \bm{\varphi}\,\dx\,\dt\\
&\quad \quad+ \int_{\T^N} \bm{\varphi}\,\int_0^{\tau} \D M^1_{E}(t) \,\dx+ \int_{0}^{\tau} \int_{\T^N}  \nabla_x \bm{\varphi}: \D\mu_m
\end{aligned}
\end{equation}
holds for all $\tau \in [0,T)$, $\p$-a.s., and for all $\bm{\varphi} \in C^{\infty}(\T^N;\mathbb{R}^N)$, where $\mu_m\in \mathcal{M}([0,T]\times{\T^N}; \mathbb{R}^N\times \mathbb{R}^N)$, $\p$-a.s., is a tensor valued measure; both $\mu_c, \mu_m$ are called {concentration measures};

\item there exists a real valued square integrable continuous martingale $M^2_{E}$, such that the following energy inequality
\begin{align}\label{measure valued energy inequality}
	\begin{aligned}
		&\mathcal{E}(t+)+ \int_{\tau}^{t} \int_{\T^N} \mathbb{S} (\nabla_x  {\langle\nu_{(r,x,\omega)},{\bf u}\rangle}) : \nabla_x {\langle\nu_{(r,x,\omega)},{\bf u}\rangle} \,\dx\,\D r\\&\quad\leq\mathcal{E}(\tau-)+\frac{1}{2} \int_\tau^{t} \bigg(\int_{\T^N} \sum_{k = 1}^{\infty} \left\langle \nu_{(r,x,\omega)};\varrho^{-1}| \mathbf{G}_k (\varrho, {\bf m}) |^2 \right\rangle \dx\bigg)\, {\rm d}r\\&\hspace{5mm}\qquad + \frac12\int_\tau^{t} \int_{\T^N} \D\mu_e + \int_{\tau}^t  \D M^2_{E}
	\end{aligned}
\end{align}

holds $\p$-a.s. for all $0\leq \tau<t$ in $(0,T)$ with \[\mathcal{E}(\tau-)\coloneqq\lim\limits_{\epsilon\to0}\frac{1}{\epsilon}\int_{\tau-\epsilon}^\tau\left(  \int_{\T^N} \left\langle \nu_{(r, x,\omega)}; \frac{1}{2} \varrho |\textbf{u}|^2 +P(\varrho) \right \rangle \dx+\mathcal{D}(r)\right)\D r,\] \[\mathcal{E}(t+)\coloneqq\lim\limits_{\epsilon\to0}\frac{1}{\epsilon}\int_{t}^{t+\epsilon}\left(  \int_{\T^N} \left\langle \nu_{(r, x,\omega)}; \frac{1}{2} \varrho |\textbf{u}|^2 +P(\varrho) \right \rangle \dx+\mathcal{D}(r)\right)\D r,\] and \[\mathcal{E}(0+)\coloneqq\int_{\T^N}\left( \frac{|\vm(0, \cdot)|^2 }{2 \vr(0, \cdot)} +  P(\vr(0, \cdot)) \right) \dx,\]

where $\mu_e\in L^\infty_{w^*}([0,T];\mathcal M_b(\T^N))$, $\p$-a.s., and $\expe{\esssup_{0\le t \le T}\mathcal{D}(t)}<\infty$, $\mathcal{D}\geq 0$, $\p$-a.s., is called {dissipation defect} of the total energy. Here the real valued martingale $M^2_E$ satisfies
\[
\E\bigg[ \sup_{t \in [0,T]} |M^2_E|^p \bigg] \leq c(p) \left( 1 + \expe{\int_{\T^N}\left( \frac{|\varrho \textbf{u}(0, \cdot)|^2 }{2 \vr(0, \cdot)} +  P(\vr(0, \cdot)) \right) \dx}^p
\right),
\]
for any $1 \leq p < \infty$; In addition, the following version of ``Poincar\'e's inequality'' holds $\p$-a.s., for a.e. $\tau \in (0,T)$
\begin{align}
\label{poincare}
\int_{0}^{\tau} \int_{\T^N} \big< \nu_{(t,x,\omega)}; |{\bf u} - {\langle\nu_{(t, x,\omega)};{\bf u}\rangle}|^2 \big> \,\dx\,\dt \le c_p\,\mathcal{D}(\tau),
\end{align}

\item there exists a constant $c>0$ such that, $\p$-a.s.
\begin{equation} \label{bound on defect measures}
\int_{0}^{\tau} \int_{\T^N} \D|\mu_c| + \int_{0}^{\tau} \int_{\T^N} \D|\mu_m| + \int_{0}^{\tau} \int_{\T^N} \D|\mu_e| \leq c \int_{0}^{\tau} \mathcal{D}(t) \,\dt,
\end{equation}	
for every $\tau \in (0,T)$.
\item Given any ``smooth'' stochastic process
$$
\mathrm{d}f  = D^d_tf\,\mathrm{d}t  + \mathbb{D}^s_tf\,\mathrm{d}W,
$$
the cross variation is given by
\begin{align*}
\left \langle f,  M^1_{E}  \right \rangle
= \int_0^T \bigg(\sum_{k = 1}^{\infty} \left\langle \nu_{(t,x,\omega)};\mathbb{D}^s_tf \,\mathbf{G}_k (\varrho,\bf{m})\right\rangle\bigg)\,\dt.
\end{align*}
\end{enumerate}
\end{Definition}	
	
\begin{Remark}
 Note that, a standard Lebesgue point argument applied to \eqref{measure valued energy inequality} gives us the usual energy inequality for a.e. $0\leq \tau<t$ in $(0,T)$ : \begin{align}\label{usual energy}
 	\begin{aligned}
 		&  \int_{\T^N} \left\langle \nu_{(t, x,\omega)}; \frac{1}{2} \varrho |\textbf{u}|^2 +P(\varrho) \right \rangle \dx + \mathcal{D}(t)+ \int_{\tau}^{t} \int_{\T^N} \mathbb{S} (\nabla_x  {\langle\nu_{(r,x,\omega)},\bf u\rangle}) : \nabla_x {\langle\nu_{(r,x,\omega)},\bf u\rangle}\dx\D r\\
 		&\quad\leq  \int_{\T^N}  \left\langle \nu_{(\tau,x,\omega)}; \frac{1}{2} \varrho |\textbf{u}|^2 +P(\varrho)\right \rangle  \dx +\mathcal{D}(\tau)+ \frac12\int_\tau^{t} \int_{\T^N} \D\mu_e + \int_{\tau}^t  \D M^2_{E}\\
 		&\qquad+ \frac{1}{2} \int_\tau^{t}\int_{\T^N} \sum_{k = 1}^{\infty} \left\langle \nu_{(r,x,\omega)};\varrho^{-1}| \mathbf{G}_k (\varrho, {\bf m}) |^2 \right\rangle \dx\,{\rm d}r\hspace{5mm}\mathbb{P} \text{-a.s.}
 	\end{aligned}			
 \end{align}
 However, for technical reasons, as pointed out in Section \ref{weakstrong uniqueness}, we require the energy inequality to hold for all $\tau,t\in(0,T)$. This is achieved via the argument described at the end of the Subsection~\ref{compactness and a.s. representations}.
 
\end{Remark}

\subsection{Main Results}
First, we have the following theorem regarding the existence of dissipative measure valued martingale solutions.
\begin{Theorem}[Existence] \label{existence of dmv martingale solution}
	Suppose that $\gamma >{\max\{1,\frac N2\}}$ and $\vc{G}_k$ be continuously differentiable satisfying the conditions \eqref{noise condition 1}--\eqref{noise condition 2}. Let $(\varrho_{\delta}, \varrho_{\delta}{\bf u}_{\delta})$ be a family of finite energy weak martingale solutions in the sense of Definition~{\normalfont\ref{finite energy weak martingale solution}} to the stochastic compressible Navier--Stokes system \eqref{P1NS}--\eqref{P2NS}. Let the corresponding initial data be $\left(\varrho_0, \varrho_0{\bf u}_{0}\right)$ and the initial law be $\Lambda$, given on the space $L^\gamma (\T^N) \times L^{\frac{2 \gamma}{\gamma + 1}}(\T^N)$, be independent of $\delta$ satisfying
	\begin{align*}
		\Lambda \Bigg\{(\varrho,\mathbf{q})\in L^\gamma(\T^N)\times L^{2\gamma/\gamma+1}(\T^N) \, :\, 0<M_1\leq \varrho\leq M_2,\, \mathbf{q}\vert_{\{\varrho=0\}}=0   \Bigg\}=1,  
	\end{align*}
	with constants $0<M_1<M_2$. Furthermore, assume that the following moment estimate
	\begin{align}\label{initialnse}
		\int_{L^\gamma_x\times  L^{2\gamma/\gamma+1}_x} \bigg\Vert\frac{1}{2}\frac{\vert  \mathbf{q} \vert^2}{\varrho}+P(\varrho)  \bigg\Vert^p_{L^1_x}\, \mathrm{d}\Lambda (\varrho,\mathbf{q})<\infty
	\end{align}
	holds for all $1\leq p<\infty$. 
	
	Then the family ${\lbrace \varrho_{\delta}, \varrho_{\delta}{\bf u}_{\delta}  \rbrace}_{\delta>0}$
	generates, as $\delta \to 0$, a Young measure ${\lbrace \nu_{(t,x,\omega)} \rbrace}_{t\in [0,T]; x\in \T^N; \omega \in \Omega}$ which is a dissipative measure valued martingale solution to
	the stochastic compressible Navier--Stokes system \eqref{stochastic ns density equation}--\eqref{stochastic ns momentum equation}, in the sense of Definition~{\normalfont\ref{dissipative measure valued martingale solution}}, with initial data $\nu_{(0,x,\omega)}=\delta_{\varrho_{0}(x,\omega), \varrho_{0}{\bf u}_0(x,\omega)}$.
\end{Theorem}

Proceeding ahead, we establish the following weak(measure valued)--strong uniqueness principle.
\begin{Theorem}[Weak--Strong Uniqueness]\label{weakstrong uniqueness for dmv solution nse}
	Let $\big[ \big(\Omega,\mathfrak{F}, (\mathfrak{F}_{t})_{t\geq0},\mathbb{P} \big); \nu_{(t,x,\omega)}, W \big]$ be a dissipative measure valued martingale solution in the sense of Definition {\normalfont\ref{dissipative measure valued martingale solution}} to the system \eqref{stochastic ns density equation}--\eqref{stochastic ns momentum equation}. On the same stochastic basis $\big(\Omega,\mathfrak{F}, (\mathfrak{F}_{t})_{t\geq0},\mathbb{P} \big)$, let us also consider the unique maximal strong path-wise solution to the Navier--Stokes system \eqref{stochastic ns density equation}--\eqref{stochastic ns momentum equation} in the sense of Definition {\normalfont\ref{def maximal strong pathwise solution for nse}}, the quadruplet $(\bar{\varrho},\bar{{\bf u}},(\mathfrak{t}_R)_{R\in\mn},\mathfrak{t})$ which is driven by the same cylindrical Wiener process $W$ with the random initial data 
	\begin{equation*} 
		\nu_{(0,x,\omega)}= \delta_{\bar{\varrho}(0,x,\omega),(\bar{\varrho}\bar{\bf u})(0,x,\omega)}, \quad \mbox{for {\normalfont a.e.} }(x,\omega)\in \T^N \times \Omega.
	\end{equation*}
	Then $\p$-{\normalfont a.s.} $\mathcal{D}(t)=0$, for {\normalfont a.e.} $t\in[0,T]$, and 
	\begin{equation*}
		\nu_{(t \wedge \mathfrak{t},x,\omega)}= \delta_{\bar{\varrho}(t \wedge \mathfrak{t},x,\omega), (\bar{\varrho}\bar{\bf u})(t \wedge \mathfrak{t},x,\omega)}, \quad \mbox{for {\normalfont a.e.} }(t,x,\omega)\in (0,T)\times \T^N \times \Omega.
	\end{equation*}
\end{Theorem}

\section[Proof of Theorem 2.10]{Proof of Theorem \ref{existence of dmv martingale solution}}\label{main proof of existence}
Thanks to the previous result {(e.g. Theorem~\ref{thm:BrFeHobook})}, the existence of weak martingale solution $$\big[ \big(\Omega_{\delta},\mathfrak{F}_{\delta}, (\mathfrak{F}_{{\delta},t})_{t\geq0},\mathbb{P}_{\delta} \big);\varrho_{\delta},\mathbf{u}_{\delta}, W_{\delta} \big]$$ of the stochastic Navier--Stokes system {\eqref{P1NS}--\eqref{P2NS}}, is well established; {see for instance \cite[Chapter~4]{BrFeHobook}.} 
The functions $\bfu_\delta$ and $\varrho_\delta$ satisfy the energy inequality, i.e., for any $1\leq p<\infty$ and {for $N\geq2$;} $\eta=\lambda+\frac{(N-2)\nu}N$ we have
\begin{align}
	&\stred\bigg[\sup_{0\leq t\leq T}\int_{\T^N}\Big(\frac{1}{2}\varrho_\delta|\bu_\delta|^2+P_{\delta}(\varrho_\delta)\Big)\,\dif x+ \int_0^T\int_{\T^N}\left(\nu|\nabla_x\bu_\delta|^2+\eta|\diver_x\bu_\delta|^2\right)\,\dif x\,\dif t\bigg]^p\nonumber\\
	&\leq \,c\,\bigg(1+\int_{L^\gamma_x\times L^\frac{2\gamma}{\gamma+1}_x}\bigg\|\frac{1}{2}\frac{|{\bf m}|^2}{\varrho}+P_{\delta}(\varrho)\bigg\|_{L^1_x}^p\,\dif \Lambda(\varrho,{\bf m})\bigg)\leq c(p,T).\label{eq:apriorivarepsilon-ch3}
\end{align}
This gives the following uniform bounds :
\begin{align}
	\bfu_\delta&\in L^{p}(\Omega;L^2(0,T;W^{1,2}( \mathbb T^N))),\label{apv-ch3}\\
	\sqrt{ \varrho_\delta} \bfu_\delta&\in L^{p}(\Omega;L^\infty(0,T;L^2( \mathbb T^N))),\label{aprhov-ch3}\\
	\varrho_\delta&\in L^{p}(\Omega;L^\infty(0,T;L^\gamma( \mathbb T^N))),\label{aprho-ch3}\\   \varrho_\delta\bfu_\delta&\in L^{p}(\Omega;L^\infty(0,T;L^\frac{2\gamma}{\gamma+1}( \mathbb T^N))),\\
	\varrho_\delta\bfu_\delta\otimes\bfu_\delta&\in L^p(\Omega;L^2(0,T;L^\frac{6\gamma}{4\gamma+3}(\T^N))).
\end{align}
Besides, the density $\varrho_\delta$ enjoys extra regularity (depending on $\delta$)
\begin{align}\label{est:nablarho}
	\delta \varrho_\delta&\in L^{p}(\Omega;L^\infty(0,T;L^\Gamma( \mathbb T^N ))).
\end{align}

{It's worth noting that, when $N=1$, the integrand $\nu|\nabla_x\bu_\delta|^2+\eta|\diver_x\bu_\delta|^2$ in the energy inequality \eqref{eq:apriorivarepsilon-ch3} boils down to $\lambda|\partial_x\bu_\delta|^2$, which makes no difference to the uniform bounds derived above.}

\subsection{Compactness and almost sure Representations}
\label{compactness and a.s. representations}

Let us define the path space $\mathcal{X}$ be the product of the following spaces :
\begin{align*}
	\mathcal{X}_\varrho&=C_w([0,T];L^\gamma(\T^N)),&\mathcal{X}_\bu&=\big(L^2(0,T;W^{1,2}(\T^N)),w\big)\\
	\mathcal{X}_{\varrho\bu}&=C_w([0,T];L^\frac{2\gamma}{\gamma+1}(\T^N)),&\mathcal{X}_W&=C([0,T];\mathfrak{U}_0)\\
	\mathcal{X}_E&= \big(L^{\infty}(0,T; \mathcal{M}_b(\T^N)), w^* \big),
	&\mathcal{X}_{P} &= \big(L^{\infty}(0,T; \mathcal{M}_b(\T^N)), w^* \big)\\
	\mathcal{X}_{C}, \mathcal{X}_{D} &= \big(L^{\infty}(0,T; \mathcal{M}_b(\T^N)), w^* \big), &\mathcal{X}_{\nu} &= \big(L^{\infty}((0,T)\times \T^N; \mathcal{P}(\R^{N+1})), w^* \big)\\
	\mathcal{X}_M&=\big(W^{\alpha, p}(0,T; W^{-m,2}(\T^N)), w \big), & \mathcal{X}_N&=C([0,T]; \R)
\end{align*}
Let us denote by $\mu_{\varrho_\delta}$, $\mu_{\bu_\delta}$, $\mu_{\varrho_\delta\bu_\delta}$, and $\mu_{W_\delta}$ respectively, the law of $\varrho_\delta$, $\bu_\delta$, $\varrho_\delta\bu_\delta$ and $W_\delta$ on the corresponding path spaces. Moreover, let $\mu_{M_\delta}$, and $\mu_{N_\delta}$ denote the laws of the martingales $M_\delta(t)=\int_0^t {\tn{G}}(\varrho_\delta, \vc{m}_\delta )\,\Dif W_\delta$, and $N_\delta(t)=\sum_{k\geq1}\int_0^t\intTorN{ \vu_{\delta} \cdot{\mathbf{G}_k}(\varrho_\delta, \vc{m}_\delta ) } \,\Dif W_{\delta,k}$. Furthermore, $\mu_{E_\delta}$, $\mu_{\nu_\delta}$, $\mu_{C_\delta}$, and $\mu_{P_\delta}$ denote the law of $E_\delta = \frac{1}{2}\varrho_\delta|\bu_\delta|^2+ P(\varrho_{\delta}) +\delta\varrho_\delta\log\varrho_\delta+ \frac{\delta}{\Gamma-1}\varrho^\Gamma_\delta$, 
$\nu_\delta = \delta_{[\varrho_\delta, \vc{m}_\delta]}$, $C_\delta= \frac{\vc{m}_\delta\otimes \vc{m}_\delta}{\varrho_\delta}$, and 
$P_\delta= p(\varrho_{\delta}) +\delta\left(\varrho_\delta+\varrho^\Gamma_\delta\right)$. $\mu_{D_\delta}$ denotes the law of $D_{\delta}= \varrho_{\delta}^{-1}| \mathbf{G}_k (\varrho_{\delta}, \vc{m}_\delta) |^2$ and their joint law on $\mathcal{X}$ is denoted by $\mu^\delta$.

To proceed, it is necessary to establish the tightness of $\{\mu^\delta;\,\delta\in(0,1)\}$. To this end, we observe that tightness of $\mu_{W_\delta}$ is immediate. We show tightness of the other variables.
\begin{Proposition}\label{tightness of velocity}
	The set $\{\mu_{\bu_\delta};\,\delta\in(0,1)\}$ is tight on $\mathcal{X}_\bu$.
\end{Proposition}

\begin{proof} 
	The proof follows directly from \eqref{apv-ch3}. Indeed, for any $R>0$ the set
	$$B_R=\big\{\bu\in L^2(0,T;W^{1,2}(\T^N));\, \|\bu\|_{L^2(0,T;W^{1,2}(\T^N))}\leq R\big\}$$
	is relatively compact in $\mathcal{X}_\bu$ and
	\begin{equation*}
		\begin{split}
			\mu_{\bu_\delta}(B_R^c)=\mathbb{P}\big(\|\bu_\delta\|_{L^2(0,T;W^{1,2}(\T^N))}\geq R\big)\leq\frac{1}{R}\stred\|\bu_\delta\|_{L^2(0,T;W^{1,2}(\T^N))}\leq \frac{c}{R},
		\end{split}
	\end{equation*}
	which yields the claim.
\end{proof}

\begin{Proposition}\label{tightness of density}
	The set $\{\mu_{\varrho_\delta};\,\delta\in(0,1)\}$ is tight on $\mathcal{X}_\varrho$. 
\end{Proposition}

\begin{proof}
	Due to \eqref{aprhov-ch3} and \eqref{aprho-ch3} we obtain that
	\begin{equation}\label{estrhou}
		\{\varrho_\delta\bfu_\delta\}\quad\text{is bounded in}\quad L^p(\Omega;L^\infty(0,T;L^{\frac{2\gamma}{\gamma+1}}(\T^N))),
	\end{equation}
	hence $\{\diver_x(\varrho_\delta\bu_\delta)\}$ is bounded in $L^{p}(\Omega;L^\infty(0,T;W^{-1,\frac{2\gamma}{\gamma+1}}(\T^N)))$. As a consequence,
	$$\stred\|\varrho_\delta\|^{p}_{C^{0,1}\left([0,T];W^{-1,\frac{2\gamma}{\gamma+1}}(\T^N)\right)}\leq c,$$
	due the continuity equation.
	Now, the required tightness in $C_w([0,T];L^\gamma(\T^N))$ follows by a similar reasoning as in Proposition \ref{tightness of velocity} together with the compact embedding (see \cite[Corollary B.2]{ond})
	$$L^\infty(0,T;L^\gamma(\T^N))\cap C^{0,1}([0,T];W^{-1,\frac{2\gamma}{\gamma+1}}(\T^N)){\compemb} C_w([0,T];L^\gamma(\T^N)).$$\end{proof}

\begin{Proposition}
	The set $\{\mu_{\varrho_\delta\bu_\delta};\,\delta\in(0,1)\}$ is tight on $\mathcal{X}_{\varrho\bu}$.
\end{Proposition}

\begin{proof}
	We first decompose $\varrho_\delta\bu_\delta$ into two parts, namely, $\varrho_\delta\bu_\delta(t)=Y^\delta(t)+Z^\delta(t)$, where
	\begin{equation*}
		\begin{split}
			Y^\delta(t)&={\bf m}_\delta(0)-\int_0^t\big[\diver_x(\varrho_\delta\bu_\delta\otimes\bu_\delta)-\nu\Delta\bu_\delta-\eta\nabla_x\diver_x\bfu_\delta
			+a\nabla_x \varrho_\delta^\gamma+\delta\nabla_x \varrho_\delta^\Gamma+\delta\nabla_x \varrho_\delta\big]\dif s,\\
			Z^\delta(t)&=\int_0^t\,\mathbb{G}(\varrho_\delta, {\bf m}_\delta) \,\dif W_\delta(s).
		\end{split}
	\end{equation*}
	Thanks to uniform bounds, we obtain H\"older continuity of $Y^\delta$, namely, there exist $\vartheta>0$ and $m>N/2$ such that $\stred\big\|Y^\delta\|_{C^\vartheta([0,T];W^{-m,2}(\T^N))}\leq c.$

	Concerning the stochastic integral, we apply BDG inequality to obtain {for some $\theta>2$,}
	\begin{equation*}
		\begin{split}
			\stred&\,\bigg\|\int_s^t\mathbb{G}(\varrho_\delta,\varrho_\delta\bu_\delta)\,\dif W_\delta\bigg\|^\theta_{W^{-m,2}(\T^N)}\leq c\,\stred\bigg(\int_s^t\sum_{k\geq1}\big\| \mathbf{G}_k(\varrho_\delta,\varrho_\delta\bu_\delta)\big\|_{W^
				{-m,2}(\T^N)}^2\,\dif r\bigg)^{\theta/2}\\
			&\leq c\,\stred\bigg(\int_s^t\sum_{k\geq1}\big\|\mathbf{G}_k(\varrho_\delta,\varrho_\delta\bu_\delta)\big\|_{L^
				{1}}^2\,\dif r\bigg)^\frac{\theta}{2}\leq c\,\stred\bigg(\int_s^t\int_{\T^N}(\varrho_\delta+\varrho_\delta|\bu_\delta|^2+\varrho_\delta^\gamma)\,\dif x\,\dif r\bigg)^{\theta/2}\\
			&\leq c|t-s|^{\theta/2}\Big(1+\stred\sup_{0\leq t\leq T}\|\sqrt\varrho_\delta\bu_\delta\|_{L^{2}}^{\theta}+\stred\sup_{0\leq t\leq T}\|\varrho_\delta\|_{L^\gamma}^{\theta\gamma/2}\Big)\leq c|t-s|^{\theta/2},
		\end{split}
	\end{equation*}
	and the Kolmogorov continuity criterion in Lemma~\ref{lemma01} applies and we get $\mathbb{E}\|Z^\delta\|_{C^\vartheta([0,T];W^{-m,2}(\T^N))}\le c$.

	Next, let us define the sets
	\begin{equation*}
		\begin{split}
			B_{R}=&\big\{h\in L^\infty(0,T;L^\frac{2\gamma}{\gamma+1}(\T^N));\,\|h\|_{L^\infty(0,T;L^\frac{2\gamma}{\gamma+1}(\T^N))}\leq R\big\},\\
			C_{R}=&\big\{h\in C^\vartheta([0,T];W^{-m,2}(\T^N));\,\|h\|_{C^\vartheta([0,T];W^{-m,2}(\T^N))}\leq R\big\},
		\end{split}
	\end{equation*}
	and
	$$K_R=B_{R}\cap C_R.$$
	Then it can be shown that $K_R$ is relatively compact in $\mathcal{X}_{\varrho\bu}$. The proof is based on the Arzel\`a-Ascoli theorem and follows closely the lines of the proof of \cite[Corollary B.2]{ond}.
	As a consequence, we obtain
	\begin{equation*}
		\begin{split}
			&\mu_{\varrho_\delta\bu_\delta}\big(K_R^c)=\mathbb{P}\big([\varrho_\delta\bu_\delta\notin B_R]\cup[Y^\delta+Z^\delta\notin C_R]\big)\\
			&\quad\leq\mathbb{P}\left(\|\varrho_\delta\bu_\delta\|_{L^\infty\left(0,T;L^\frac{2\gamma}{\gamma+1}(\T^N)\right)}>R\right)+\mathbb{P}\Big(\|Y^\delta\|_{C^\vartheta([0,T];W^{-m,2}(\T^N))}>\frac R2\Big)\\&\qquad+\mathbb{P}\Big(\|Z^\delta\|_{C^\vartheta([0,T];W^{-m,2}(\T^N))}>\frac R2\Big)\leq \frac{c}{R}.
		\end{split}
	\end{equation*}
	A suitable choice of $R$ completes the proof.
\end{proof}

\begin{Proposition}\label{tightness1}
	The set $\{\mu_{E_{\delta}}, \mu_{C_{\delta}}, \mu_{D_{\delta}},  \mu_{P_{\delta}};\,\delta\in(0,1)\}$ is tight on $\mathcal{X}_{E} \times \mathcal{X}_{C} \times \mathcal{X}_{D}  \times \mathcal{X}_{P}$.
\end{Proposition}

\begin{proof}
	These follow immediately from the {\it a priori} bounds and using the fact that all bounded sets in $L^{\infty}(0,T;\mathcal{M}_b(\T^N))$ are relatively compact with respect to the weak$^*$ topology.
\end{proof}

\begin{Proposition}\label{tightness2}
	The set $\{\mu_{\nu_{\delta}};\,\delta\in(0,1)\}$ is tight on $\mathcal{X}_{\nu}$.
\end{Proposition}

\begin{proof}
	The aim here is to apply the compactness criterion in the space $\big(L^{\infty}((0,T)\times \T^N; \mathcal{P}(\R^{N+1})), w^* \big)$.

	Define the set
	\begin{align*}
		B_R= \left\lbrace \nu \in L^{\infty}_{w^*}\left((0,T)\times \T^N; \mathcal{P}\left(\R^{N+1}\right)\right)\Bigg| 
		\int\limits_{(0,T)\times\T^N} \int\limits_{\R^{N+1}} \Big(|\xi_1|^{\gamma} + |\xi_2|^{\frac{2\gamma}{\gamma+1}} \Big)\dif\nu_{(t,x)}\dx\dt \le R    \right\rbrace,
	\end{align*}
	which is relatively compact in $\big(L^{\infty}((0,T)\times \T^N; \mathcal{P}(\R^{N+1})), w^* \big)$. Note that
	\begin{align*}
		\mathcal{L}[\nu_{\delta}](B^c_R)&=
		\p\left( \int_0^T \int_{\T^N} \int_{\R^{N+1}} \Big(|\xi_1|^{\gamma} + |\xi_2|^{\frac{2\gamma}{\gamma+1}} \Big) \,\D\nu_{(t,x)}(\xi)\,\dx\,\dt > R  \right) \\
		&= \p\left(\int_0^T \int_{\T^N} \left(|\varrho_\delta|^{\gamma} + |\vc{m}_{\delta}|^{\frac{2\gamma}{\gamma+1}} \right) \,\dx\,\dt >R \right)\\
		&\le \frac1R \E\left[\|\varrho_\delta \|_{\gamma}^{\gamma} + \| \vc{m}_{\delta}\|_{\frac{2\gamma}{\gamma+1}}^{\frac{2\gamma}{\gamma+1}}\right] \le \frac cR.
	\end{align*}
	The proof is complete.\end{proof}

\begin{Proposition}
	The set $\{\mu_{M_{\delta}};\,\delta\in(0,1)\}$ is tight on $\mathcal{X}_{M}$.
\end{Proposition}

\begin{proof}
	{Using \cite[Lemma~2.1]{FlandoliGatarek}, we directly obtain for $p\ge2$ and $\alpha<\frac12$, $$M_{\delta}= \int_0^t\,\mathbb{G}(\varrho_\delta,\varrho_\delta\bu_\delta) \,\dif W_\delta(s) \in L^p\big(\Omega; W^{\alpha, p}(0,T; W^{-m,2}(\T^N))\big),$$} and hence tightness follows from weak compactness.
\end{proof}

\begin{Proposition}
	The set $\{\mu_{N_{\delta}};\,\delta\in(0,1)\}$ is tight on $\mathcal{X}_{N}$.
\end{Proposition}

\begin{proof}
	First note that, for each $\delta$, $N_\delta(t)=\sum_{k\geq1}\int_0^t\intTorN{ \vu_{\delta} \cdot {\mathbf{G}_k}(\varrho_\delta, \vc{m}_\delta ) } \,\Dif W_{\delta,k}$ is a square integrable martingale and {by the BDG inequality, one observes} for $r>2$,
	\begin{align*}
		\E\Big[ \Big|\int_s^t \int_{\T^N}\vu_{\delta} \cdot {\mathbf{G}_k}(\varrho_\delta, \vc{m}_\delta )\dx\,\Dif W_{\delta,k}\Big|^r  \Big] &\le c\,\E\Big[ \int_s^t \sum_{k=1}^{\infty} \Big|\int_{\T^N}\vu_{\delta} \cdot {\mathbf{G}_k}(\varrho_\delta, \vc{m}_\delta)\dx\Big|^2\Dif\tau  \Big]^{r/2} \\
		& \le c |t-s|^{r/2}\, \Big(1 + \E \sup_{0\le t \le T} \| \sqrt{\varrho_\delta} u_{\delta}\|^r_{2} \Big)\\
		&\le c|t-s|^{r/2},
	\end{align*}
	and the Kolmogorov continuity criterion (cf. Lemma~\ref{lemma01}) applies. This in particular implies, for some {$0<\alpha<\frac{1}{2}-\frac{1}{r}$,}
	$$
	\sum_{k\geq1}\int_0^t\intTorN{ \vu_{\delta} \cdot {\mathbf{G}_k}(\varrho_\delta, \vc{m}_\delta ) } \,\Dif W_{\delta,k} \in L^r(\Omega; C^{\alpha}(0,T; \R)).
	$$
	Therefore, the tightness of laws follows from the compact embedding of $C^{\alpha}$ into $C^0$.
\end{proof}

\begin{Corollary}
	The set $\{\mu^\delta;\,\delta\in(0,1)\}$ is tight on $\mathcal{X}$. 
\end{Corollary}

Now we have all in hand to apply the Jakubowski-Skorokhod representation theorem. It yields the following.

\begin{Proposition}\label{prop:skorokhod1-ch3}
	There exists a subsequence $\mu^\delta$, a probability space $(\tilde\Omega,\tilde{\mathfrak{F}},\tilde{\mathbb{P}})$ with $\mathcal{X}$-valued Borel measurable random variables $(\tilde\varrho_\delta,\tilde\bu_\delta,\tilde{\mathbf{m}}_\delta,\tilde W_\delta,\cdots)$, and $(\tilde\varrho,\tilde\bu,\tilde{\mathbf{m}},\tilde W,\cdots)$ such that
	\begin{enumerate}
		\item[\normalfont1.] the law of $(\tilde\varrho_\delta,\tilde\bu_\delta,\tilde{\bf m}_\delta,\tilde W_\delta, \cdots)$ is given by $\mu^\delta$, $\delta\in(0,1)$,
		\item[\normalfont2.] the law of $(\tilde\varrho,\tilde\bu,\tilde{\bf m},\tilde W, \cdots)$, denoted by $\mu$, is a Radon measure,
		\item[\normalfont3.] $(\tilde\varrho_\delta,\tilde\bu_\delta,\tilde{\bf m}_\delta,\tilde W_\delta,\cdots)$ converges $\,\tilde{\mathbb{P}}$-almost surely to $(\tilde\varrho,\tilde{\bu},\tilde{\bf m},\tilde{W},\cdots)$ in the topology of $\mathcal{X}$, i.e.,
		\begin{align*}
			&\tilde\varrho_\delta \rightarrow \tilde\varrho \,\, \text{in}\, \,C_w([0,T]; L^{\gamma}(\T^N)),\\
			&\tilde{\vc{m}}_\delta \rightarrow \tilde{\vc m} \,\, \text{in}\, \,C_w([0,T]; L^{\frac{2\gamma}{\gamma+1}}(\T^N)),\\
			&\tilde{\vc{u}}_\delta \rightarrow \tilde{\vc u} \,\, \text{weak in}\, \,L^2( 0,T; W^{1,2}(\T^N)),\\
			&\tilde W_\delta \rightarrow \tilde W \,\, \text{in}\, \,C([0,T]; \mathfrak{U}_0),\\
			& \tilde E_\delta \rightarrow \tilde E \,\, \text{weak$^*$ in}\, \, L^{\infty}(0,T; \mathcal{M}_b(\T^N)),\\
			& \tilde C_\delta \rightarrow \tilde C \,\, \text{weak$^*$ in}\, \, L^{\infty}(0,T; \mathcal{M}_b(\T^N)),\\
			& \tilde D_\delta \rightarrow \tilde D \,\, \text{weak$^*$ in}\, \, L^{\infty}(0,T; \mathcal{M}_b(\T^N)),\\
			& \tilde P_\delta \rightarrow \tilde P \,\, \text{weak$^*$ in}\, \, L^{\infty}(0,T; \mathcal{M}_b(\T^N)),\\
			& \tilde \nu_\delta \rightarrow \tilde \nu \,\, \text{weak$^*$ in}\, \, L^{\infty}((0,T)\times \T^N; \mathcal{P}(\R^{N+1})),\\
			& \tilde M_\delta \rightarrow \tilde M \,\, \text{weak in}\, \, W^{\alpha, p}(0,T; W^{-m,2}(\T^N)),\\
			&\tilde N_\delta \rightarrow \tilde N \,\, \text{in}\, \, C([0,T]; \R),
		\end{align*}
		\item[\normalfont4.] For any Carath\'{e}odory function $H=H(t,x,\varrho,\vc m)$, where $(t,x)\in (0,T)\times \T^N$ and $(\varrho,\vc m) \in \R\times\R^{N}$, satisfying for some $p, q$ the growth condition
		\begin{align*}
			|H(t,x,\varrho,\vc m)| \le 1 + |\varrho|^{p} + |\vc m|^q,
		\end{align*}
		uniformly in $(t,x)$. Then we have
		$$
		H(\tilde\varrho_\delta, \tilde{\vc{m}}_\delta) \rightarrow \overline{H(\tilde\varrho, \tilde{\vc m})}\,\, \text{in}\,\, L^r((0,T)\times\T^N),\,\, \text{for all}\,\, 1<r\le\frac{\gamma}{p}\wedge \frac{2\gamma}{q(\gamma+1)}
		$$
		as $\delta \to 0$, $\tilde\p${\normalfont-a.s.}
	\end{enumerate}
\end{Proposition}
\begin{proof}
	Proof of the items $(1)$, $(2)$, and $(3)$ directly follow from Jakubowski--Skorokhod representation theorem (cf. Theorem~\ref{JakubowskiSkorohod}). For the proof of the item $(4)$, we refer to the Lemma~\ref{younglemma}.
\end{proof}

We remark that the energy inequality \eqref{eq:apriorivarepsilon-ch3} continues to hold on the new probability space. In other words, all the {\it a priori} estimates \eqref{apv-ch3}--\eqref{est:nablarho} also hold for the new random variables.

Let $(\tilde{\mathfrak{F}}_t^\delta)$ and $(\tilde{\mathfrak{F}}_t)$, respectively, be the $\tilde{\mathbb{P}}$-augmented canonical filtration of the process $(\tilde\varrho_\delta,\tilde{\bf m}_\delta,\tilde{W}_\delta)$ and $(\tilde\varrho,\tilde{\bf m},\tilde{W})$, respectively, that is
\begin{equation*}
	\begin{split}
		\tilde{\mathfrak{F}}_t^\delta&=\sigma\big(\sigma\big(\bfr_t\tilde\varrho_\delta,\,\bfr_t\tilde{\bf m}_\delta,\,\bfr_t \tilde{W}_\delta\big)\cup\big\{N\in\tilde{\mathfrak{F}};\;\tilde{\mathbb{P}}(N)=0\big\}\big),\quad t\in[0,T],\\
		\tilde{\mathfrak{F}}_t&=\sigma\big(\sigma\big(\bfr_t\tilde\varrho, \,\bfr_t\tilde{\bf m},\,\bfr_t\tilde{W}\big)\cup\big\{N\in\tilde{\mathfrak{F}};\;\tilde{\mathbb{P}}(N)=0\big\}\big),\quad t\in[0,T].
	\end{split}
\end{equation*}

\begin{Proposition}
	For every $\delta\in(0,1)$, $\big((\tilde{\Omega}_\delta,\tilde{\mathfrak{F}}_\delta,(\tilde{\mathfrak{F}}_{\delta,t}),\tilde{\mathbb{P}}_\delta),\tilde\varrho_\delta,\tilde{\bu}_\delta,\tilde{W}_\delta\big)$ is a finite energy weak martingale solution to \eqref{P1NS}--\eqref{P2NS} with the initial law $\Lambda_\delta$. 
\end{Proposition}

\begin{proof}
	To see this, we follow the monograph by Breit and Hofmanová \cite[Proposition 4.11]{BrHo}, and define accordingly, for all $t\in[0,T]$ and $\bfvarphi\in C^\infty(\T^N;\R^N)$ the functionals
	\begin{equation*}
		\begin{split}
			\tilde M_\delta(\tilde\varrho_\delta,\tilde\bfv_\delta,\tilde\bfq_\delta)_t&=\big\langle \tilde\bfq_\delta(t),\bfvarphi\big\rangle-\big\langle \tilde\bfq_\delta(0),\bfvarphi\big\rangle-\int_0^t\big\langle \tilde\bfq_\delta\otimes \tilde\bfv_\delta,\nabla_x\bfvarphi\big\rangle\,\dif r+\nu\int_0^t\big\langle\nabla_x \tilde\bfv_\delta,\nabla_x\bfvarphi\big\rangle\,\dif r\\
			&\hspace{2mm}+\eta\int_0^t\big\langle\diver_x \tilde\bfv_\delta,\diver_x\bfvarphi\big\rangle\,\dif r-a\int_0^t\big\langle \tilde\varrho_\delta^\gamma,\diver_x\bfvarphi\big\rangle\,\dif r-\delta\int_0^t\big\langle \tilde\varrho_\delta^\Gamma+\tilde\varrho_\delta,\diver_x\bfvarphi\big\rangle\,\dif r,\\
			N_\delta(\tilde\varrho_\delta,\tilde\bfq_\delta)_t&=\sum_{k\geq1}\int_0^t\big\langle \mathbf{G}_k(\tilde\varrho_\delta, \tilde\bfq_\delta ),\bfvarphi\big\rangle^2\,\dif r,\\
			N_{k,\delta}(\tilde\varrho_\delta,\tilde\bfq_\delta)_t&=\int_0^t\big\langle \mathbf{G}_k(\tilde\varrho_\delta,\tilde\bfq_\delta),\bfvarphi\big\rangle\,\dif r,
		\end{split}
	\end{equation*}
	Below $\tilde M_\delta(\tilde\varrho_\delta,\tilde\bfv_\delta,\tilde\bfq_\delta)_{s,t}$ denotes the increment $\tilde M_\delta(\tilde\varrho_\delta,\tilde\bfv_\delta,\tilde\bfq_\delta)_t-\tilde M_\delta(\tilde\varrho_\delta,\tilde\bfv_\delta,\tilde\bfq_\delta)_s$ and similarly for other processes and then we deduce that
	\begin{equation}
		\begin{split}
			&\tilde{\stred}\,h\big(\bfr_s\tilde\varrho_\delta, \bfr_s\tilde{\bu}_\delta,\bfr_s\tilde{W}_\delta\big)\big[M_{\delta}(\tilde\varrho_\delta,\tilde\bu_\delta,\tilde\varrho_\delta\tilde\bfu_\delta)_{s,t}\big]=0,\\
			&\tilde{\stred}\,h\big(\bfr_s\tilde\varrho_\delta, \bfr_s\tilde{\bu}_\delta,\bfr_s\tilde{W}_\delta\big)\bigg[[M_{\delta}(\tilde\varrho_\delta,\tilde\bu_\delta,\tilde\varrho_\delta\tilde\bfu_\delta)^2]_{s,t}-N_\delta(\tilde\varrho_\delta,\tilde\varrho_\delta\tilde\bfu_\delta)_{s,t}\bigg]=0,\\
			&\tilde{\stred}\,h\big(\bfr_s\tilde\varrho_\delta, \bfr_s\tilde{\bu}_\delta,\bfr_s\tilde{W}_\delta\big)\bigg[[M_{\delta}(\tilde\varrho_\delta,\tilde\bu_\delta,\tilde\varrho_\delta\tilde\bfu_\delta)\tilde{W}_{\delta,k}]_{s,t}-N_{k,\delta}(\tilde\varrho_\delta,\tilde\varrho_\delta\tilde\bfu_\delta)_{s,t}\bigg]=0,
		\end{split}
	\end{equation}
	which essentially prove the result.\end{proof}

 The above proposition implies that the new random variables satisfy:
\begin{itemize}
	\item for all $\phi\in C^\infty(\T^N)$ and $\bm{\varphi}\in C^\infty(\T^N;\R^N)$ we have
	\begin{equation}\label{eq:energyt-ch3}
		\begin{aligned}
			&\langle \tilde\varrho_{\delta}(t); \phi\rangle = \langle\tilde\varrho_{\delta}(0) ; \phi\rangle + \int_0^{t}\langle \tilde {\bf m}_{\delta};\nabla_x \phi\rangle\mathrm{d}s, 
			\\
			&\langle \tilde {\bf m}_{\delta}(t), \bm\varphi\rangle = \langle \tilde {\bf m}_{\delta}(0);\bm{\varphi}\rangle + \int_0^{t} \bigg\langle \frac{\tilde {\bf m}_{\delta}\otimes\tilde {\bf m}_{\delta}}{\tilde \varrho_{\delta}}; \nabla_x  \bm{\varphi} \bigg \rangle\,\mathrm{d}s
			- \int_0^{t} \langle\mathbb{S}(\nabla_x \tilde{\mathbf{u}}_{\delta}); \nabla_x  \bm{\varphi}\rangle\mathrm{d}s 
			\\&\hspace{20mm}
			+\int_0^{t}\langle  p_{\delta}(\tilde\varrho_{\delta}), \Div\bm{\varphi}\rangle\,\mathrm{d}s
			+\int_0^{t}\langle\mathbb{G}(\tilde \varrho_{\delta},\tilde {\bf m}_{\delta}), \bm{\varphi}\rangle\,\mathrm{d}\tilde{W}_\delta,
		\end{aligned}
	\end{equation}
	$\tilde{\mathbb{P}}$-a.s. for all $t\in[0,T]$,
	\item the energy inequality\index{energy inequality}
	\begin{align}\label{EI2t''-ch3}
		\begin{aligned}
			&
			-\int_0^T \partial_t\psi\int_{\T^N} \bigg[ \frac{1}{2} \frac{ | \tilde {\bf m}_{\delta} |^2 }{\tilde\varrho_{\delta}} + P_{\delta}(\tilde\varrho_{\delta}) \bigg] \dx\,\dt
			+ \int_0^T \psi\int_{\T^N}  \mathbb{S} (\nabla_x  \tilde {\bf u}_{\delta}): \nabla_x  \tilde {\bf u}_{\delta} \dx\,\dt\\
			& \leq \psi(0)\int_{\T^N} \bigg[ \frac{1}{2}\frac{|\tilde {\bf m}_{\delta}(0)|^2}{\tilde\varrho_{\delta}(0)}  + P_{\delta}(\tilde\varrho_{\delta})(0) \bigg] \dx
			+\sum_{k=1}^\infty\int_0^T\psi\int_{\T^N}\mathbf{G}_k (\tilde \varrho_{\delta}, \tilde {\bf m}_{\delta})\cdot \tilde{\bf u}_{\delta}\dx\,{\rm d} \tilde{W}_{\delta,k}\\
			&\qquad + \frac{1}{2}\sum_{k = 1}^{\infty}  \int_0^T\psi
			\int_{\T^N} \tilde \varrho_{\delta}^{-1}| \mathbf{G}_k (\tilde \varrho_{\delta}, \tilde {\bf m}_{\delta}) |^2 \dx\, {\rm d}t,
		\end{aligned}
	\end{align}
	$\tilde{\mathbb P}$-a.s. holds for any $\psi\in C_c^\infty([0,T))$, with $\psi\geq0$. 
\end{itemize}

\subsubsection{Passage to Limits}
 Next, our aim is to pass to the limit in $\delta$ in the approximate equations \eqref{eq:energyt-ch3} and the energy inequality \eqref{EI2t''-ch3}. To do this, we make use of the Young measure theory. First note that, thanks to previous results, $\tilde \p$-a.s.,
\begin{align*}
	&\tilde \varrho_\delta \rightarrow \langle \tilde\nu_{(t,x,\omega)}; \tilde \varrho \rangle\,\,\text{in}\,\, L^{\gamma}((0,T)\times\T^N),\\
	&\tilde {\vc m}_\delta \rightarrow \langle \tilde\nu_{(t,x,\omega)}; \tilde {\vc m} \rangle\,\,\text{in}\,\, L^{\frac{2\gamma}{\gamma+1}}((0,T)\times\T^N).
\end{align*}
In order to pass to the limit in the nonlinear terms, we can introduce new measures and follow the deterministic techniques. Let 
\begin{align*}
	&\tilde \mu_{C}= \tilde C -\left\langle \tilde\nu_{(t,x,\omega)}; \frac{\tilde {\bf m}\otimes \tilde {\bf m}}{\tilde\varrho} \right\rangle \dx\dt,\\&\tilde \mu_{P}= \tilde P -\langle \tilde\nu_{(t,x,\omega)};p(\tilde\varrho)\rangle \dx\dt,\\&\tilde \mu_{E}= \tilde E- \left\langle \tilde\nu_{(t,x,\omega)}; \frac{1}{2} \frac{|\tilde {\bf m}|^2}{\tilde\varrho} +P(\tilde\varrho) \right \rangle \dx\dt,\\&\tilde \mu_{D}= \tilde D -\left\langle \tilde \nu_{(t,x,\omega)}; \sum_{k \geq 1} \frac{ |{\bf G}_k (\tilde \varrho, \tilde {\bf m}) |^2 }{\tilde \varrho}\right\rangle \dx\dt.
\end{align*}
For the extra term coming from the pressure perturbation, we argue as follows: Note that $\tilde\p$-a.s. 
$$
\delta\left(\tilde{\varrho}_\delta\log\tilde{\varrho}_\delta+\frac{\tilde{\varrho}_{\delta}^{\Gamma}}{\Gamma-1}\right) \rightharpoonup \xi, \,\, \mbox{weak$^*$ in}\,\, L^{\infty}((0,T); \mathcal{M}_b^+(\T^N)).
$$
Observe that to obtain the above convergence, we need to put this variable in the list of variables before applying the Skorokhod theorem. Therefore, while passing to the limit in the energy term, a new term $\xi(\tau)$ will appear due to the perturbation of the pressure. We shall absorb this term in the energy dissipation defect measure $\mathcal{D}(\tau)$ and then we are able to pass to the limit in the dissipation term on the L.H.S of the energy inequality by using the regularity of $\tilde{\bf u}_{\delta}$. In fact, 
\begin{align*}
	&\mathbb{S} (\nabla_x  \tilde{\bf u}_{\delta}) : \nabla_x {\tilde{\bf u}_{\delta}} = \nu |\nabla_x {\tilde{\bf u}_{\delta}}|^2 + \eta |\Div {\tilde{\bf u}_{\delta}}|^2\\&\hspace{28mm}\overset{\ast}{\rightharpoonup}\nu |\nabla_x \big<\nu_{(t,x,\omega)}; {\bf u} \big>|^2 + \eta \big( \mathrm{tr} |\nabla_x \big<\nu_{(t,x,\omega)}; {\bf u} \big>|\big)^2+\sigma_\infty(t),
\end{align*}where $\sigma_\infty(t)$ is a non-negative measure on $\T^N$, which is again absorbed in $\mathcal{D}(t)$.

Again, new terms coming from the momentum equation can be dealt with the same strategy. Indeed, here also the defect measure is dominated by the energy dissipation defect.

Regarding the martingale $\tilde M_{\delta}$ that appears in the momentum equation, we obtain {from Lemma~\ref{FlandoliGatarek}}, the compact embedding {(by choosing $p=4,\alpha=\frac13$ so that $\alpha p=\frac{4}{3}>1$)} \[W^{\alpha, p}(0,T; W^{-m,2}(\T^N))\compemb C\left([0,T]; W^{-m-1,2}(\T^N)\right).\] This implies, for each $t$, $\tilde M_{\delta}(t) \rightarrow \tilde M(t)$, $\tilde\p$-a.s. in the topology of $W^{-m-1,2}(\T^N)$. Therefore, we can pass to the limit in $\delta$ to get the desired result (cf. \cite[Proposition~3.10]{mku}). For the other martingale term appearing in the energy inequality, we conclude that, after a change of the probability space, a new sequence $\{ \tilde N_{\delta} \}_{\delta=\frac1k;k\geq 1}$ having the same law as the original sequence $\{ N_{\delta} \}_{\delta=\frac1k;k\geq 1}$, converges to some $\tilde N$ a.s. in $\mathcal{X}_N$ (cf. \cite[Proposition~3.9]{mku}). Since the space of continuous square integrable martingales is closed, we deduce that the limits $\tilde{M},\tilde N$ are also square integrable martingales. Besides, it follows from the equality of the joint laws that the energy inequality is also satisfied on the new probability space.\vspace{2mm}

Now collecting all the above information, we infer that, \begin{align*}
	&\int_{\T^N}\left\langle\tilde{\nu}_{(\omega,t,x)};\varrho\right\rangle\phi(x)\,\dx=\int_{\T^N}\left\langle\tilde{\nu}_{(\omega,0,x)};\varrho\right\rangle\phi(x)\,\dx+\int_0^t	\int_{\T^N}\left\langle\tilde{\nu}_{(\omega,s,x)};\vm\right\rangle\cdot\nabla_x\phi\,\dx\,{\D s}
\end{align*}
holds $\tilde{\mathbb{P}}$-a.s. for all $t\in[0,T)$, and for all $\phi\in C^\infty(\T^N)$. Notice here, due to the regularity of $\tilde{\varrho}_\delta,\tilde{\vm}_\delta$, we have no concentration measure, i.e., $\tilde{\mu}_c=0$. 

Moreover, \begin{align*}
	&\int_{\T^N}\left\langle\tilde{\nu}_{(\omega,t,x)};\vm\right\rangle\cdot{\bm\varphi}(x)\,\dx-\int_{\T^N}\left\langle\tilde{\nu}_{(\omega,0,x)};\vm\right\rangle\cdot{\bm\varphi}(x)\,\dx\\&\hspace{5mm}=-\int_0^t\int_{\T^N}\mathbb{S}\left(\nabla_x\left\langle\tilde{\nu}_{(\omega,s,x)};\vu\right\rangle\right):\nabla_x\bm\varphi \,\dx\,\D s\\
	&\hspace{8mm}+\int_0^t\int_{\T^N}\left[\left\langle\tilde{\nu}_{(\omega,s,x)};\frac{\vm\otimes\vm}{\varrho}\right\rangle:\nabla_x\bm{\varphi}+\left\langle\tilde{\nu}_{(\omega,s,x)};p(\varrho)\right\rangle\Div\bm{\varphi}\right]\dx\,\D s\\&\hspace{8mm}+\int_0^t\int_{\T^N}\nabla_x\bm{\varphi}:\D\left(\tilde{\mu}_C+\tilde{\mu}_P\mathbb{I}\right)+\int_{\T^N}\bm{\varphi}(x)\int_0^t\D \tilde{M}(s)\,\dx
\end{align*} holds\,\,$\tilde{\mathbb{P}}$-a.s. for all $t\in[0,T)$, and for all $\bm{\varphi}\in C^\infty\left(\T^N;\R^N\right)$, where $\tilde{\mu}_m\coloneqq\tilde{\mu}_C+\tilde{\mu}_P\mathbb{I}\in L^\infty_{w^*}\left([0,T];\mathcal{M}_b(\T^N)\right)$, $\tilde{\mathbb{P}}$-a.s. is a tensor valued measure, and thus we conclude \eqref{measure valued density equation} and \eqref{measure valued momentum equation}. To see \eqref{measure valued energy inequality}, we proceed as follows. First note, after letting $\delta\to0$ in \eqref{EI2t''-ch3}, we get the following energy inequality in the new probability space :\begin{align}\label{energyfinal}
	\begin{aligned}
		&-\int_0^T \partial_t\psi\left[\int_{\T^N} \left\langle\tilde{\nu}_{(\omega,t,x)};\frac{ | {\bf m}|^2 }{2\varrho} + P(\varrho)\right\rangle \dx+\tilde{\mathcal{D}}(t)\right]\dt\\
		&\hspace{10mm}+ \int_0^T \psi\int_{\T^N} \mathbb{S} (\nabla_x  {\langle\nu_{(\omega,t,x)},{\bf u}\rangle}) : \nabla_x \langle\nu_{(\omega,t,x)},{\bf u}\rangle\dx\,\dt\\
		&\hspace{3mm}\leq \psi(0)\int_{\T^N} \left\langle\tilde{\nu}_{(\omega,0,x)};\frac{|{\bf m}|^2}{2\varrho}  + P(\varrho) \right\rangle \dx+\frac{1}{2}\int_0^T\psi\int_{\T^N}\D\tilde\mu_D+\int_0^T\psi\,\D\tilde{N}\\
		&\hspace{10mm}+ \frac{1}{2} \int_0^T\psi
		\int_{\T^N} \sum_{k = 1}^{\infty} \left\langle\tilde{\nu}_{(\omega,t,x)};\varrho^{-1}| \mathbf{G}_k (\varrho, \tilde {\bf m}) |^2\right\rangle \dx\, {\rm d}t.
	\end{aligned}
\end{align}
Here, according to our previous notations $0\leq\tilde{\mathcal{D}}(t)\coloneqq\tilde{\mu}_E(t)(\T^N),\,\tilde{\mu}_D=\tilde{\mu}_e,\,\tilde{N}=\tilde{M}^2_E$ and $\tilde{M}=\tilde{M}^1_E$. Now we choose a specific function $\psi_r:[0,T]\mapsto\R$ as follows: Fix any $t_0,t$ such that $0<t_0<t<T$, and for any $r>0$ such that $0<t_0-r<t+r<T$, we define $\psi_r$ to be the function that is linear on $[t_0-r,t_0]\cup[t,t+r]$ and satisfies \[\psi_r(s)=\begin{cases}
	0\hspace{3mm}\text{if}\hspace{1mm}s\in[0,t_0-r]\cup[t+r,T],\\
	1\hspace{3mm}\text{if}\hspace{1mm}s\in[t_0,t].
\end{cases}\]
Then, via a standard regularization argument, we see that $\psi_r$ is an admissible test function in \eqref{energyfinal}. Inserting this $\psi_r$ into \eqref{energyfinal} gives \begin{align*}
	&\frac{1}{r}\int_t^{t+r}\left[\int_{\T^N} \left\langle\tilde{\nu}_{(\omega,s,x)};\frac{ | {\bf m}|^2 }{2\varrho} + P(\varrho)\right\rangle \dx+\tilde{\mathcal{D}}(s)\right]\D s\\
	&\hspace{10mm}+ \int_{t_0-r}^{t+r} \int_{\T^N} \mathbb{S} (\nabla_x  {\langle\nu_{(\omega,s,x)},{\bf u}\rangle}) : \nabla_x \langle\nu_{(\omega,s,x)},{\bf u}\rangle\dx\,\D s\\
	&\hspace{3mm}\leq \frac{1}{r}\int_{t_0-r}^{t_0}\left[\int_{\T^N} \left\langle\tilde{\nu}_{(\omega,s,x)};\frac{ | {\bf m}|^2 }{2\varrho} + P(\varrho)\right\rangle \dx+\tilde{\mathcal{D}}(s)\right]\D s\\
	&\hspace{10mm}+ \frac{1}{2} \int_{t_0-r}^{t+r}
	\int_{\T^N} \sum_{k = 1}^{\infty} \left\langle\tilde{\nu}_{(\omega,s,x)};\varrho^{-1}| \mathbf{G}_k (\varrho, \tilde {\bf m}) |^2\right\rangle \dx\, {\rm d}s+\frac{1}{2}\int_{t_0-r}^{t+r}\int_{\T^N}\D\tilde\mu_D+\int_{t_0-r}^{t+r}\,\D\tilde{N}.
\end{align*}
Now, taking limit as $r\to0^+$, we conclude that \eqref{measure valued energy inequality} holds. Note that for $t_0=0$, we take $\psi_r$, that assumes $1$ in $[0,t]$, linear on $[t,t+r]$ and zero otherwise and apply the same argument as above.

\begin{Proposition}
	Given any ``smooth'' stochastic process
	$$
	\mathrm{d}f  = D^d_tf\,\mathrm{d}t  + \mathbb{D}^s_tf\,\mathrm{d}W,
	$$
	the cross variation is given by
	\begin{align*}
		\left \langle f,  \tilde M \right \rangle
		= \int_0^T \bigg( \sum_{k=1}^\infty\Big\langle \tilde \nu_{(t,x,\omega)};\mathbb{D}^s_tf \,\mathbf{G}_k (\tilde \varrho,\tilde {\bf m})\Big\rangle \bigg)\,\dt.
	\end{align*}
\end{Proposition}

\begin{proof}
	Note that we have
	\begin{align*}
		\left \langle f,  \tilde M_{\delta} \right \rangle = \sum_{k = 1}^{\infty} \int_0^T \Big\langle \mathbf{G}_k (\tilde \varrho_{\delta},\tilde {\bf m}_{\delta}), \mathbb{D}^s_tf  \Big \rangle \,\dt.
	\end{align*}
	We may apply part $(4)$ of Proposition~\ref{prop:skorokhod1-ch3} to the composition $\mathbf{G}_k (\tilde \varrho_{\delta},\tilde {\bf m}_{\delta})$, $k \in \mathbb{N}$. This gives
	$$
	\mathbf{G}_k (\tilde \varrho_{\delta},\tilde {\bf m}_{\delta}) \rightarrow
	\Big \langle \tilde \nu_{(t,x,\omega)} ; \mathbf{G}_k (\tilde \varrho,\tilde {\bf m}) \Big \rangle \,\, \mbox{in} \,\,L^q((0,T)\times \T^N),
	$$
	$\tilde\p$-a.s., for some $q>1$. Moreover, for $m>N/2$, we have by Sobolev embedding 
	\begin{align*}
		\E \int_0^T \| \mathbb{G}(\tilde{\varrho}_{\delta}, \tilde{\bf m}_{\delta}) \|^2_{L_2(\mathfrak{U}; W^{-m,2})}\,\dt \le \E \int_0^T (\tilde{\varrho}_{\delta})_{\T^N} \int_{\T^N} (\tilde{\varrho}_{\delta} +\tilde {\varrho}_{\delta} |\tilde{\bf u}_{\delta}|^2)\,\dx\,\dt \le c.
	\end{align*}
	This implies that
	$$
	\mathbf{G}_k (\tilde \varrho_\delta,\tilde {\bf m}_\delta) \rightharpoonup \Big \langle \tilde \nu_{(t,x,\omega)} ; \mathbf{G}_k (\tilde \varrho,\tilde {\bf m}) \Big \rangle \,\, \mbox{weakly in} \,\,L^2((0,T); W^{-m,2}(\T^N)),
	$$ 
	$\tilde\p$-a.s. Now Vitali lemma gives the required convergence (cf. \cite[Lemma~3.11]{mku}).
\end{proof}

\begin{Proposition}
	The defect measure $0\leq\tilde{\mathcal{D}}(t)=\tilde \mu_{E}(t)(\T^N)$ dominates defect measures $\tilde{\mu}_C,\,\tilde{\mu}_P,$ and $\tilde \mu_{D}$ in the sense of Lemma~{\normalfont\ref{lemma001}}. Specifically, \eqref{bound on defect measures} holds.
\end{Proposition}

\begin{proof} Following deterministic argument, we can conclude that $\tilde{\mu}_{E}$ dominates the defect measures $\tilde{\mu}_C,\,\tilde{\mu}_P$.
	Next we observe, by virtue of the hypotheses (\ref{noise condition 1}), (\ref{noise condition 2}), the function
	\[
	[\vr, \vc{m}] \mapsto \sum_{k \geq 1} \frac{ |{\bf G}_k (\varrho, {\bf m}) |^2 }{\varrho} \ \mbox{is continuous},
	\]
	and
	\[
	\sum_{k \geq 1} \frac{ |{\bf G}_k (\varrho, {\bf m}) |^2 }{\varrho} \leq c \left(  \varrho  + \frac{| {\bf m} |^2}{\vr}  \right)
	\]
	is sub-linear in $\varrho$ and $|{\bf m}|^2/\vr$ and therefore is dominated by the total energy
	\[
	\left( \frac{|{\bf m}|^2}{2\vr}  + P(\vr) \right) + 1.
	\]
	Hence the required result follows.
\end{proof}


\section{Weak--Strong Uniqueness Principle}
\label{weakstrong uniqueness}

We start by recalling the following auxiliary result from \cite[Lemma~3.1]{BrFeHo2015A}, which is essential in our proof.

\begin{Lemma} \label{ito-product lemma}

Let $w$ be a stochastic process on $\StoB$ such that for some $\lambda>0$,
\[
w \in C_{\rm weak}([0,T]; W^{-\lambda,2}(\tor)) \cap L^\infty (0,T; L^1(\tor)) \quad \text{$\mathbb{P}$-{\normalfont a.s.}},
\]
\begin{equation} \label{hh1}
 \E\bigg[ \sup_{t \in [0,T]} \| w \|^p_{L^1(\tor)}\bigg]  < \infty \ \mbox{for all}\ 1 \leq p < \infty,
\end{equation}
\begin{equation} \label{rel1}
\Dif w= \DD w\, \dt + \DS w\,\Dif W.
\end{equation}
Here $\DD w, \DS w$ are progressively measurable with
\begin{align} \label{hh2}
\begin{aligned}
\DD w\in &L^p(\Omega;L^1(0,T;W^{-\lambda,q}(\tor)),\quad \DS w\in L^2(\Omega;L^2(0,T;L_2(\mathfrak U;W^{-m,2}(\tor)))),\\
&\sum_{k\geq 1}\int_0^T\|\DS w(e_k)\|^2_1\in L^p(\Omega)\quad 1\leq p<\infty,
\end{aligned}
\end{align}
for some $q>1$ and some $m\in\N$.

Let $r$ be a stochastic process on $\StoB$ satisfying
\[
r \in C([0,T]; W^{\lambda,q'} \cap C (\tor)) \quad \text{$\mathbb{P}$-{\normalfont a.s.}},
\]
\begin{equation} \label{hh3}
\E\bigg[ \sup_{t \in [0,T]} \| r \|_{W^{\lambda,q'} \cap C (\tor)}^p \bigg] < \infty,\ 1 \leq p < \infty,
\end{equation}
\begin{equation} \label{rel2}
\Dif r = \DD r \,\dt+ \DS r \, \Dif W.
\end{equation}
Here $q'=\frac{q}{q-1}\mbox{ and }\DD r, \DS r$ are progressively measurable with
\begin{align} \label{hh4}
\begin{aligned}
\DD r\in L^p&(\Omega;L^1(0,T;W^{\lambda,q'}\cap C(\tor)),\quad \DS r\in L^2(\Omega;L^2(0,T;L_2(\mathfrak U;W^{-m,2}(\tor)))),\\
&\sum_{k\geq1} \int_0^T\|\DS r(e_k)\|^2_{W^{\lambda,q'}\cap C(\tor)}\dt\in L^p(\Omega),\quad 1\leq p<\infty.
\end{aligned}
\end{align}
Let $Q$ be $\lfloor\lambda+2\rfloor$-continuously differentiable function satisfying
\begin{equation} \label{hh5}
\E\bigg[\sup_{t \in [0,T]} \| Q^{(j)} (r) \|_{W^{\lambda,q'} \cap C (\tor)}^p \bigg]< \infty, \quad j = 0,1,2,\quad 1 \leq p < \infty.
\end{equation}

Then
\begin{equation} \label{result}
\begin{split}
&\Dif \left( \intTorN{ w Q(r) } \right)\\
&= \bigg( \intTorN{ \Big[  w  \Big( Q'(r) \DD r   + \frac{1}{2}\sum_{k\geq1} Q''(r)
\left| \DS r (e_k)\right|^2  \Big) \Big] }  +  \left< Q(r) , \DD w \right> \bigg) {\rm d}t
\\
&\qquad+ \bigg(  \sum_{k\geq1}\intTorN{  \DS w(e_k)\,\DS r(e_k)  } \bigg) {\rm d}t
+ {\rm d}\tn{M},
\end{split}
\end{equation}
where
\begin{equation} \label{result1}
\tn{M} = \sum_{k\geq1}\int_0^t\intTorN{ \Big[  w  Q'(r) \DS r(e_k)  + Q(r) \DS w(e_k)  \Big] }\,\Dif W_k.
\end{equation}
\end{Lemma}

\qed

\begin{Remark}

The result stated in Lemma \ref{ito-product lemma} is not optimal with respect to the regularity properties of the processes
$r$ and $w$. As a matter of fact, we could regularize both $r$ and $w$ in the above proof to conclude that \eqref{result} holds as long as all the expressions in \eqref{result} and \eqref{result1} are well-defined.

\end{Remark}

\subsection{Relative Energy Inequality}\label{subsection relative energy inequality}
We introduce the following {relative energy functional} in the context of measure valued solutions for the compressible Navier--Stokes system :
\[
\begin{split}
	\mathcal{E}_{\mathrm{mv}}^1 \left( \varrho,{\bf m} \ \Big| \ r, {\bf U} \right)(t)
	\coloneqq&
	\int_{\T^N}{ \Bigg\langle \nu_{(t,x,\omega)}; \frac{1}{2} \frac{ |{\bf m}|^2}{\vr} + P(\varrho) \Bigg\rangle}\,\D x+\mathcal{D}(t) - \int_{\T^N}{ \big \langle \nu_{(t,x,\omega)}; {\bf m} \big \rangle\cdot {\bf U}}\,\D x
	\\&+ \frac{1}{2}  \int_{\T^N}{ \big \langle \nu_{(t,x,\omega)}; \vr \big \rangle \, |{\bf U}|^2 }\,\D x - \int_{\T^N}{ \big \langle \nu_{(t,x,\omega)}; \vr \big \rangle \, P'(r) }\,\D x\\& + \int_{\T^N}{ \left[ P'(r) r - P(r) \right]\D x }.
\end{split}
\]
We remark that the above functional is defined for all $t\in[0,T]\backslash N$, where the set $N$, may depend on $\omega$, has Lebesgue measure zero. To define the functional for all $t\in N$, we consider 
\[
\begin{split}
	\mathcal{E}_{\mathrm{mv}}^2 \left( \varrho,{\bf m} \ \Big| \ r, {\bf U} \right)(t )
	\coloneqq&
	\lim\limits_{\tau\to0}\frac{1}{\tau}\int_t^{t+\tau}\left[\int_{\T^N}{ \Bigg\langle \nu_{(t,x,\omega)}; \frac{1}{2} \frac{ |{\bf m}|^2}{\vr} + P(\varrho) \Bigg\rangle}\,\D x+\mathcal{D}(s)\right]\D s\\& \hspace{1cm}- \int_{\T^N}{ \big \langle \nu_{(t,x,\omega)}; {\bf m} \big \rangle\cdot {\bf U}}\,\D x+ \frac{1}{2}  \int_{\T^N}{ \big \langle \nu_{(t,x,\omega)}; \vr \big \rangle \, |{\bf U}|^2 }\,\D x\\&\hspace{1cm} - \int_{\T^N}{ \big \langle \nu_{(t,x,\omega)}; \vr \big \rangle \, P'(r) }\,\D x+ \int_{\T^N}{ \left[ P'(r) r - P(r) \right]\D x }.
\end{split}
\]

Using above, we define the relative energy functional as follows
\begin{equation}
	\label{refnse}
	\begin{aligned}
		\mathcal{E}_{\mathrm{mv}} \left( \varrho,{\bf m} \ \Big\vert \ r, {\bf U} \right)(t) \coloneqq
		\begin{cases}
			\mathcal{E}^1_{\mathrm{mv}} \left( \varrho,{\bf m} \ \vert \ r, {\bf U} \right)(t), \, &\text{if} \,\,t \in [0,T] \setminus N,\\
			\mathcal{E}^2_{\mathrm{mv}} \left( \varrho,{\bf m} \ \vert \ r, {\bf U} \right)(t), \, &\text{if} \,\,t \in N.
		\end{cases}
	\end{aligned}
\end{equation}
Note that the above relative energy functional is well defined for all $t \in [0,T]$. In what follows, with the help of the above defined relative energy, we derive the relative energy inequality \eqref{relativeEntropy}. Note that, the relative energy inequality is a tool which enables us to compare measure valued solutions with some smooth comparison functions.\vspace{3mm}

We have the following inequality:

\begin{Theorem}[Relative Energy]
Let $\big[ \big(\Omega,\mathfrak{F}, (\mathfrak{F}_{t})_{t\geq0},\mathbb{P} \big); \nu_{(t,x,\omega)}, W \big]$ be a dissipative measure valued martingale solution to the system \eqref{stochastic ns density equation}--\eqref{stochastic ns momentum equation}.
Suppose $\big(r \,,\,\mathbf{U}  \big)$ be a pair of stochastic processes which are adapted to the filtration $(\mathfrak{F}_t)_{t\geq0}$ and satisfy
\begin{equation}
\begin{aligned}
\label{operatorBB}
\mathrm{d}r  &= D^d_tr\,\mathrm{d}t  + \mathbb{D}^s_tr\,\mathrm{d}W,
\\
\mathrm{d}\mathbf{U}  &= D^d_t\mathbf{U}\,\mathrm{d}t  + \mathbb{D}^s_t\mathbf{U}\,\mathrm{d}W,
\end{aligned}
\end{equation}
with
\begin{align}\label{rU}
\begin{aligned}
r \in C([0,T]; W^{1,q}(\T^N)), \ \vc{U} \in C([0,T]; W^{1,q}(\T^N)), \ \quad\text{$\mathbb{P}$-{\normalfont a.s.}},\\
\E\bigg[\sup_{t \in [0,T] } \| r \|_{W^{1,q}(\tor)}^2\bigg]^q  + \E\bigg[ \sup_{t \in [0,T] } \| \vc{U} \|_{W^{1,q}(\tor)}^2\bigg]^q \leq c(q) \quad\forall\,\, 2 \leq q < \infty,
\end{aligned}
\end{align}
\begin{equation} \label{bound on smooth process}
0 < \underline{r} \leq r(t,x) \leq \overline{r} \quad\text{$\mathbb{P}$-{\normalfont a.s.}}
\end{equation}
Moreover, $r$, $\vc{U}$ satisfy
\begin{align}
&D^d r, D^d \vc{U}\in L^q(\T^N;L^q(0,T;W^{1,q}(\mt))),\quad
\mathbb D^s r,\mathbb D^s \vc{U}\in L^2(\T^N;L^2(0,T;L_2(\mathfrak U;L^2(\tor)))),\nonumber\\
&\bigg(\sum_{k\geq 1}|\mathbb D^s r(e_k)|^q\bigg)^\frac{1}{q},\bigg(\sum_{k\geq 1}|\mathbb D^s \vc{U}(e_k)|^q\bigg)^\frac{1}{q}\in L^q(\T^N;L^q(0,T;L^{q}(\mt))).\label{bound on smooth process 2}
\end{align}
Then the following {relative energy inequality} holds $\mathbb P$-{\normalfont a.s.} for all $t\in[0,T]$:
\begin{equation}
\begin{aligned}
\label{relativeEntropy}
&\mathcal{E}_{\mathrm{mv}} \left(\varrho,\textbf{m} \ \Big| \ r, {\bf U} \right)
(t)\\
&\quad+\int_0^t\int_{\T^N}\left(\mathbb{S}\left(\nabla_x{\langle\nu_{(s,x,\omega)},\bf u\rangle}\right)-\mathbb{S}(\nabla_x\bf U)\right):\left(\nabla_x{\langle\nu_{(s,x,\omega)},\bf u\rangle}-\nabla_x\bf U\right)\dx\,\D s \\&\leq
\mathcal{E}_{\mathrm{mv}} \left(\varrho,\textbf{m} \ \Big| \ r, {\bf U} \right)(0) +M_{RE}(t)  + \int_0^t\mathcal{R}_{\mathrm{mv}} \big(\varrho,\textbf{m} \left\vert \right. r, \mathbf{U}  \big)(s)\,\D s,
\end{aligned}
\end{equation}
 where
\begin{equation}
\begin{aligned}
\label{remainderRE}
\mathcal{R}_{\mathrm{mv}} \big(\varrho,{\bf m} \left\vert \right. r, \mathbf{U}  \big)\,\dt
&=\int_{\T^N}\mathbb{S}\left(\nabla_x{\bf U}\right):\left(\nabla_x{\bf U}-\nabla_x\left\langle\nu_{(t,x,\omega)},\frac{{\bf m}}{\varrho}\right\rangle\right)\dx\,\dt
\\&+\int_{\T^N}\langle \nu_{(t,x,\omega)}; \varrho \mathbf{U} - {\bf m} \rangle \cdot[D^d_t  \mathbf{U} + \mathbf{U}\cdot\nabla_x  \mathbf{U}] \,\dx\,\dt \\
&+ \int_{\T^N} \left\langle \nu_{(t,x,\omega)}; \frac{({\bf m}-\varrho  \mathbf{U})\otimes (\varrho  \mathbf{U}-{\bf m}) }{\varrho} \right\rangle: \nabla_x {\bf U} \,\dx\,\dt
\\
&+
\int_{\T^N}\big[(r- \big \langle \nu_{(t,x,\omega)}; \varrho \big\rangle)P''(r)D^d_tr + \nabla_x  P'(r) \cdot (r\mathbf{U}-\big \langle \nu_{(t,x,\omega)}; {\bf m} \big\rangle)  \big] \,\dx\,\dt
\\
&+
\int_{\T^N}\big[p(r)  - \big \langle \nu_{(t,x,\omega)}; p(\varrho)\big\rangle \big]\Div\mathbf{U} \,\dx\,\dt
\\
&+
\frac{1}{2}
\sum_{k\in\mathbb{N}}
\int_{\T^N} \Bigg \langle \nu_{(t,x,\omega)}; \varrho\bigg\vert \frac{\mathbf{G}_k(\varrho,{\bf m})}{\varrho}  -\mathbb{D}^s_t\mathbf{U}(e_k)  \bigg\vert^2\Bigg\rangle \,\dx\,\dt \\
&-  \int_{\T^N} \nabla_x \mathbf{U} : \D\mu_m + \frac12\int_{\T^N} \D \mu_e\\
&-\frac{1}{2}\int_{\T^N}\sum_{k \in\mathbb{N}}\langle\nu_{(t,x,\omega)},\varrho\rangle P'''(r)\lvert\mathbb{D}^s_tr\rvert^2\dx\,\dt+\frac12 \int_{\T^N}\sum_{k \in\mathbb{N}}p''(r)\lvert\mathbb{D}^s_tr\rvert^2\dx\,\dt.
\end{aligned}
\end{equation}
Here $M_{RE}$ is a real valued square integrable martingale, and the norm of this martingale depends only on the norms of $r$ and $\mathbf{U}$ in the aforementioned spaces. Moreover the pressure potential $P$ is defined as the solution of the equation $p(r)=rP'(r)-P(r)$.

\end{Theorem}

\begin{proof}
As the time evolution of the first integral is governed by the energy inequality \eqref{measure valued energy inequality}, it remains to compute the time differentials
of the remaining terms with the help of Lemma \ref{ito-product lemma}.

\medskip

 {\bf Step 1:}

\medskip

To compute $\Dif \intTorN{ \big \langle \nu_{(t,x,\omega)}; \textbf{m} \big \rangle \,\cdot {\bf U} }$ we recall that $w = \big \langle \nu_{(t,x,\omega)}; \textbf{m} \big \rangle $ satisfies hypotheses \eqref{hh1}--\eqref{hh2} with $m=1$ and some $1<q < \infty$. Applying Lemma \ref{ito-product lemma} we obtain
\begin{equation} \label{ito step 1}
\begin{split}
& \Dif \left( \intTorN{ \left \langle \nu_{(t,x,\omega)}; \textbf{m} \right\rangle \cdot {\bf U} } \right)\\
& =\intTorN{ \left[ \left \langle \nu_{(t,x,\omega)}; \textbf{m} \right \rangle \cdot \DD {\bf U}
+ \left\langle \nu_{(t,x,\omega)}; \frac{\textbf{m}\otimes \textbf{m}}{\varrho} \right\rangle: \nabla_x \textbf{U} \right] }\,{\rm d}t \\
&\quad+\intTorN{\left \langle \nu_{(t,x,\omega)}; p(\vr) \right\rangle\, \Div {\bf U} + \nabla_x{\bf U}:\mathbb{S}\left(\nabla_x{\left\langle\nu_{(t,x,\omega)},\frac{\bf m}{\varrho}\right\rangle}\right)}\dt\\
& \qquad +  \sum_{k\geq1}\intTorN{ \DS {\bf U}(e_k) \cdot \left \langle \nu_{(t,x,\omega)}; \vc{G}_k (\varrho,\textbf{m}) \right\rangle}\, {\rm d}t + \int_{\T^N}  \nabla_x {\bf U}: \D\mu_m+ \Dif M_1,
\end{split}
\end{equation}
where
\[
M_1(t) = \int_0^t \int_{\T^N} \vc{U} \, \D M^1_{E} \,\dx+ \int_0^t \intTorN{ \big \langle \nu_{(t,x,\omega)}; \textbf{m} \big \rangle \cdot \DS \vc{U} } \,\Dif W
\]
is a square integrable martingale.

\medskip

{\bf Step 2:}

\medskip

Similarly, we compute
\begin{equation} \label{ito step 2}
\begin{split}
\Dif \left( \intTorN{ \frac{1}{2} \big \langle \nu_{(t,x,\omega)}; \vr \big \rangle |\vc{U}|^2 } \right) &=
\intTorN{ \big \langle \nu_{(t,x,\omega)}; \textbf{m} \big \rangle \cdot \Grad {\bf U} \cdot {\bf U} }  \,{\rm d}t+  \intTorN{ \big \langle \nu_{(t,x,\omega)}; \vr \big \rangle {\bf U} \cdot  \DD {\bf U} } \,{\rm d}t
\\&\quad+ \frac{1}{2} \sum_{k\geq1}\intTorN{ \big \langle \nu_{(t,x,\omega)}; \vr \big \rangle |\DS {\bf U}(e_k)|^2 } \,{\rm d}t + {\rm d}M_2,
\end{split}
\end{equation}
where
\[
M_2(t) = \int_0^t \intTorN{ \big \langle \nu_{(t,x,\omega)}; \vr \big \rangle \vc{U} \cdot \DS \vc{U} } \, \Dif W.
\]

{\bf Step 3:}

\medskip

Next we compute,
\begin{equation} \label{ito step 3}
{\rm d} \left( \intTorN{ \left[ P'(r) r - P(r) \right] } \right)
= \intTorN{ p'(r) \DD r } \ {\rm d}t + \frac{1}{2}\sum_{k\geq1} \intTorN{ p''(r) |\DS r(e_k)|^2 } \ {\rm d}t + \Dif M_3,
\end{equation}
where
\[
M_3(t) = \int_0^t \intTorN{ p'(r) \DS r } \,\Dif W,
\]

{\bf Step 4:}

\medskip

and, finally
\begin{equation} \label{ito step 4}
\begin{split}
{\rm d} \left( \intTorN{ \big \langle \nu_{(t,x,\omega)}; \vr \big \rangle P'(r) } \right) =&
  \intTorN{  \Grad P'(r) \cdot \big \langle \nu_{(t,x,\omega)}; \textbf{m} \big \rangle } \ {\rm d}t+ \intTorN{ \big \langle \nu_{(t,x,\omega)}; \vr \big \rangle P''(r) \DD r } \ {\rm d}t\\ &\quad\frac{1}{2}\sum_{k\geq1} \intTorN{ \big \langle \nu_{(t,x,\omega)}; \vr \big \rangle P'''(r) |\DS r(e_k) |^2 } \ {\rm d}t
+ \Dif M_4,
\end{split}
\end{equation}
where
\[
M_4(t) = \int_0^t \intTorN{ \big \langle \nu_{(t,x,\omega)}; \vr \big \rangle P''(r) \DS r } \,\Dif W.
\]

\medskip

{\bf Step 5:}

Now we define $M_{RE}\coloneqq M^2_E-M_1+M_2+M_3-M_4$ and collect \eqref{ito step 1}--\eqref{ito step 4} and summing up the resulting expressions and adding the sum with \eqref{measure valued energy inequality}, we obtain the inequality \eqref{relativeEntropy}.
\end{proof}

\subsection[Proof of Theorem 2.11]{Proof of Theorem \ref{weakstrong uniqueness for dmv solution nse}}

As we have seen before, the {relative energy inequality} \eqref{relativeEntropy} is a consequence of the {energy inequality} \eqref{measure valued energy inequality}. Now we wish to use the above theorem to prove the weak--strong uniqueness principle:

\begin{proof}
In order to prove \ref{weakstrong uniqueness for dmv solution nse}, we start by choosing $(r,\mathbf U)=(\bar{\varrho}(\cdot\wedge \mathfrak t_R),\bar{\bfu}(\cdot\wedge \mathfrak t_R))$, where $(\bar{\varrho},\bar{\mathbf{u}},(\mathfrak{t}_R)_{R\in\mathbb{N}},\mathfrak{t})$ is the unique maximal strong path-wise solution to \eqref{stochastic ns density equation}--\eqref{stochastic ns momentum equation} which exists by Theorem \ref{existence of finite enerey weak martingale solution for nse} . Recall that the stopping time $\mathfrak t_R$ announces the blow-up and satisfies
\begin{equation*}
\sup_{t\in[0,\mathfrak{t}_R]}\|\bar{\vu}(t)\|_{W^{2,\infty}}\geq R\quad \text{on}\quad [\mathfrak{t}<T].
\end{equation*}
Moreover, $(r,\mathbf U)=(\bar{\varrho},\bar{\bfu})$ satisfies an equation of the form \eqref{operatorBB}, where
\begin{align}\label{value}
	\begin{aligned}
	&	D^d_tr  = -\Div (\bar{\varrho}\bar{\mathbf{u}}),\quad \mathbb{D}^s_tr=0,
		\\
		&D^d_t\mathbf{U}  = -\bar{\mathbf{u}}\cdot\nabla_x \bar{\mathbf{u}} -\frac{1}{\bar{\varrho}}\nabla_x  p(\bar{\varrho})+\frac{\Div\mathbb{S}(\nabla_x\bar{\bf u})}{\bar{\varrho}},
		\quad
		\mathbb{D}^s_t\mathbf{U}  = \frac{1}{\bar{\varrho}}\mathbb{G} (\bar{\varrho},\bar{\varrho}\bar{{\bf u}}).
	\end{aligned}
\end{align}
By Theorem \ref{existence of finite enerey weak martingale solution for nse} and \eqref{noise condition 1}--\eqref{noise condition 2}, it is easy to see that
\eqref{rU} and \eqref{bound on smooth process 2} are satisfied for $t\leq \mathfrak t_R$. Note in particular the lower bound for $\bar{\varrho}$ which follows from the maximum principle and also $\bar{\varrho}$ satisfies the equation of continuity on the interval $[0,\mathfrak{t}]$, so there exists deterministic constants $\underline{\varrho}_R,\overline{\varrho}_R$ such that $$0<\underline{\varrho}_R\le\bar{\varrho}(t\wedge \mathfrak{t}_R)\le\overline{\varrho}_R\text{ on }[0,T].$$
So, \eqref{relativeEntropy} holds and we can now deduce from \eqref{remainderRE} that for each $R\in\mathbb{N}$ and every $t\in[0,T]$,
\begin{equation}
\begin{aligned}
\label{relativeEntropy1}
&\mathcal{E}_{\mathrm{mv}} \left( \varrho,\textbf{m} \,\Big| \,\bar{\varrho}, \bar{\mathbf{u}} \right)  
(t \wedge \mathfrak{t}_R)\\ &\qquad+\int_0^{t\wedge\mathfrak{t}_R}\int_{\T^N}\left(\mathbb{S}\left(\nabla_x{\langle\nu_{(s,x,\omega)},\bf u\rangle}\right)-\mathbb{S}(\nabla_x \bar{{\bf u}})\right):\left(\nabla_x{\langle\nu_{(s,x,\omega)},\bf u\rangle}-\nabla_x\bar{\bf u}\right)\dx\,\D s\\
& \qquad \qquad\leq
\mathcal{E}_{\mathrm{mv}} \left( \varrho,\textbf{m} \,\Big| \,\bar{\varrho}, \bar{\mathbf{u}} \right)(0) + M_{RE}(t \wedge \mathfrak{t}_R)  + \int_0^{t \wedge \mathfrak{t}_R} \mathcal{R}_{\mathrm{mv}} \left(\varrho,{\bf m} \,\Big\vert\,\bar{\varrho}, \bar{\mathbf{u}}  \right) (s)\,\mathrm{d}s.
\end{aligned}
\end{equation} 
Now plugging \eqref{value} in \eqref{remainderRE} and after a straightforward computation, following \cite{BrFeHo2015A,FeNoSu}, we obtain
\begin{equation}
\begin{aligned}
\label{remainderRE1}
\mathcal{R}_{\mathrm{mv}} \left(\varrho,{\bf m} \,\Big\vert\,\bar{\varrho}, \bar{\mathbf{u}}  \right)(s) \,\D s
=&\int_{\T^N}\frac{\Div\mathbb{S}\left(\nabla_x\bar{{\bf u}}\right)}{\bar{\varrho}}\cdot\langle\nu_{(s,x,\omega)};(\varrho-\bar{\varrho})(\bar{\vu}-\vu)\rangle\,\dx\,\D s
\\&+\int_{\T^N} \left\langle \nu_{(s,x,\omega)}; \frac{(\textbf{m}-\varrho \bar{\mathbf{u}})\otimes (\varrho \bar{\mathbf{u}}-\textbf{m}) }{\varrho}\right\rangle : \nabla_x \bar{\mathbf{u}}\,\dx\,\D s\\
&-
\int_{\T^N}\left\langle \nu_{(s,x,\omega)};p(\varrho)-p'(\bar{\varrho})(\varrho-\bar{\varrho})-p(\bar{\varrho})\right\rangle \divv_x \bar{\mathbf{u}}\,\dx\,\D s\\
&+
\frac{1}{2}
\sum_{k\in\mathbb{N}}
\int_{\T^N} \left\langle \nu_{(s,x,\omega)}; \varrho\left\vert \frac{\mathbf{G}_k(\varrho,\varrho \mathbf{u})}{\varrho}  -\frac{\mathbf{G}_k(\bar{\varrho}, \bar{\varrho} \bar{\mathbf{u}})}{\bar{\varrho}}  \right\vert^2 \right\rangle \,\mathrm{d}x\,\D s \\
&- \int_{\T^N} \nabla_x \bar{\mathbf{u}}:\D\mu_m+\frac{1}{2}\int_{\T^N}\D{\mu_e}.
\end{aligned}
\end{equation}

Before proceeding, we set a few notations and state an elementary inequality here, that is crucial in our proof. We denote $H(\varrho,r)\coloneqq P(\varrho)-P'(r)(\varrho-r)-P(r)$ and then, it is easy to check that for $\alpha>0$, \begin{equation}\label{inequality0}
H(\varrho,r)\geq c(\alpha)\begin{cases*}
		\left\lvert\varrho-r\right\rvert^2\hspace{2mm}\text{for}\hspace{2mm}\alpha<\varrho,r<\alpha^{-1},\\
		1+\varrho^\gamma\hspace{2mm}\text{whenever}\hspace{2mm}\alpha<r<\alpha^{-1},\varrho\in(0,\infty)\backslash\left[\frac\alpha2,2\alpha\right].
	\end{cases*}
\end{equation}

Since $\bar{\textbf{u}}$ {is a strong path-wise solution}, we control the terms $|\nabla_x \bar{\mathbf{u}}|$, $|\divv_x \bar{\mathbf{u}}|$ and $|\divv_x \mathbb{S}\left(\nabla_x\bar{\mathbf{u}}\right)|$ {uniformly in $(t,x)\in[0,\mathfrak{t}_R]\times\mathbb{T}^N$} by some constant {$c(R)$}. It is also obvious that there exist constants $c_1>0$ and $c_2>\gamma-1$ such that
\begin{equation}
\begin{aligned}
&J_2\coloneqq\left| \frac{(\textbf{m}- \varrho\bar{\mathbf{u}}) \otimes (\varrho \bar{\mathbf{u}}-\textbf{m})}{\varrho} \right| \leq \frac{c_1}{2\varrho} |\textbf{m}-\varrho \bar{\mathbf{u}}|^2,\\
&J_3\coloneqq|p(\varrho) -p'(\bar{\varrho})(\varrho - \bar{\varrho})- p(\bar{\varrho})| \leq c_2 (P(\varrho) -P'(\bar{\varrho})(\varrho - \bar{\varrho})- P(\bar{\varrho})).
\end{aligned}
\end{equation}
Which implies
\begin{equation}\label{inequality2}
\begin{aligned}
\intTorN{\left\langle\nu_{(t,x,\omega)};J_2\right\rangle|\nabla_x\bar{\vu}|}+\intTorN{\left\langle\nu_{(t,x,\omega)};J_3\right\rangle|\Div \bar{\vu}|}\leq c(R)\,\mathcal{E}_{\mathrm{mv}}\left(\varrho,\bf m\,\big\vert\,\bar{\varrho},\bar{\vu}\right).
\end{aligned}
\end{equation}
We rewrite
\begin{equation}
\begin{aligned}
\label{nio0}
&J_4\coloneqq\varrho\bigg\vert \frac{\mathbf{G}_k(\varrho,\varrho \mathbf{u})}{\varrho}  -\frac{\mathbf{G}_k(\bar{\varrho}, \bar{\varrho}\bar{\mathbf{u}})}{\bar{\varrho}}  \bigg\vert^2
= \chi_{\{\varrho \leq \bar{\varrho}/2\}}\, \varrho\bigg\vert \frac{\mathbf{G}_k(\varrho,\varrho \mathbf{u})}{\varrho}  -\frac{\mathbf{G}_k(\bar{\varrho}, \bar{\varrho} \bar{\mathbf{u}})}{\bar{\varrho}}  \bigg\vert^2
\\
& \qquad + \chi_{\{ \bar{\varrho}/2 < \varrho < 2\bar{\varrho}\}}\, \varrho\bigg\vert \frac{\mathbf{G}_k(\varrho,\varrho \mathbf{u})}{\varrho}  -\frac{\mathbf{G}_k(\bar{\varrho}, \bar{\varrho} \bar{\mathbf{u}})}{\bar{\varrho}}  \bigg\vert^2 
+ \chi_{\{  \varrho \geq 2\bar{\varrho}\}}\, \varrho\bigg\vert \frac{\mathbf{G}_k(\varrho,\varrho \mathbf{u})}{\varrho}  -\frac{\mathbf{G}_k(\bar{\varrho}, \bar{\varrho} \bar{\mathbf{u}})}{\bar{\varrho}}  \bigg\vert^2
\\
&\quad\eqqcolon I_1+I_2+I_3.
\end{aligned}
\end{equation}
We now use the inequality $\varrho<1+\varrho^\gamma$ and \eqref{inequality0} to conclude that
\begin{equation}
\begin{aligned}
\label{inequality4}
I_1 &\leq 
c\, \chi_{\{\varrho \leq \bar{\varrho}/2\}}\,\left( \frac{1}{\varrho}\,\vert \mathbf{G}_k(\varrho ,\varrho {\mathbf{u}}) \vert^2  +\frac{\varrho}{\bar{\varrho}^2}\,\vert \mathbf{G}_k(\bar{\varrho},\bar{\varrho}\bar{\mathbf{u}}) \vert^2  \right) 
\\
&\leq c\, \chi_{\{\varrho \leq  \bar{\varrho}/2\}}\left( \varrho+ \varrho\vert {\mathbf{u}}\vert^2  + \varrho\vert  \bar{\mathbf{u}}\vert^2\right)\\
&\leq c\,\chi_{\{\varrho \leq  \bar{\varrho}/2\}}\left(\left(1+\varrho^\gamma\right)(1+3|\bar{{\bf u}}|^2)+2\varrho \vert {\mathbf{u}} -  \bar{\mathbf{u}}\vert^2\right) \\
& \leq c(R)\left(P(\varrho) -P'(\bar{\varrho})(\varrho - \bar{\varrho})- P(\bar{\varrho})+\varrho \vert {\mathbf{u}} -  \bar{\mathbf{u}}\vert^2\right).
\end{aligned}
\end{equation}
Therefore,
\begin{equation}
\begin{aligned}
\sum_{k\in\mathbb{N}} \int_{\T^N} \big \langle \nu_{(t,x,\omega)}; I_1\big \rangle \,\mathrm{d}x 
\leq c(R)\,
\mathcal{E}_{\mathrm{mv}} \big(\varrho ,{\mathbf{m}}\left\vert \right. \bar{\varrho}, \bar{\mathbf{u}}  \big).
\end{aligned}
\end{equation}
We can estimate $I_3$ in an exactly similar fashion. For $I_2$, we see 
\begin{equation}\label{inequality5}
	\begin{aligned}
	I_2&\leq c\chi_{\{ \bar{\varrho}/2 < \varrho < 2\bar{\varrho}\}}\, \varrho\bigg\vert \frac{\mathbf{G}_k(\varrho,\varrho \mathbf{u})}{\varrho}  -\frac{\mathbf{G}_k(\bar{\varrho}, {\varrho} {\mathbf{u}})}{\bar{\varrho}}  \bigg\vert^2\\ &\quad+c\chi_{\{ \bar{\varrho}/2 < \varrho < 2\bar{\varrho}\}}\, \varrho\bigg\vert \frac{\mathbf{G}_k(\bar{\varrho},\varrho \mathbf{u})}{\bar{\varrho}}  -\frac{\mathbf{G}_k(\bar{\varrho}, \bar{\varrho} \bar{\mathbf{u}})}{\bar{\varrho}}  \bigg\vert^2\\
	&\leq c(R)\chi_{\{ \bar{\varrho}/2 < \varrho < 2\bar{\varrho}\}}\left(\lvert\varrho-\bar{\varrho}\rvert^2\left(1+|\varrho{\bf u}|^2\right)+|\varrho{\bf u}-\bar{\varrho}\bar{\vu}|^2\right)\\
	&\leq c(R)\chi_{\{ \bar{\varrho}/2 < \varrho < 2\bar{\varrho}\}}\left(|\varrho-\bar{\varrho}|^2(1+|\bar{\vu}|^2)+\varrho|\vu-\bar{\vu}|^2\right)\\
	&\leq c(R)\left(P(\varrho)-P'(\bar{\varrho})(\varrho-\bar{\varrho})-P(\bar{\varrho})+\varrho|\vu-\bar{\vu}|^2\right).
\end{aligned}
\end{equation}
Therefore,
\begin{equation}
	\begin{aligned}
		\sum_{k\in\mathbb{N}} \int_{\T^N} \big \langle \nu_{(t,x,\omega)}; I_2\big \rangle \,\mathrm{d}x 
		\leq c(R)\,
		\mathcal{E}_{\mathrm{mv}} \big(\varrho ,{\mathbf{m}}\left\vert \right. \bar{\varrho}, \bar{\mathbf{u}}  \big).
	\end{aligned}
\end{equation}
Finally, we can conclude from \eqref{nio0} that
\begin{equation}
\begin{aligned}
\label{inequality3}
\frac{1}{2}
\sum_{k\in\mathbb{N}}
\int_{\T^N} \Bigg \langle \nu_{(t,x,\omega)}; \varrho\bigg\vert \frac{\mathbf{G}_k(\varrho,\varrho \mathbf{u})}{\varrho}  -\frac{\mathbf{G}_k(\bar{\varrho}, \bar{\varrho} \bar{\mathbf{u}})}{\bar{\varrho}}  \bigg\vert^2 \Bigg \rangle \,\mathrm{d}x 
\leq c(R)\,
\mathcal{E}_{\mathrm{mv}} \big(\varrho ,{\mathbf{m}}\left\vert \right. \bar{\varrho}, \bar{\mathbf{u}}  \big).
\end{aligned}
\end{equation}
Using \eqref{bound on defect measures}, we get \begin{equation}\label{inequality6}
	\begin{aligned}
		\int_0^{t\wedge\mathfrak{t}_R} \int_{\T^N} |\nabla_x \bar{\mathbf{u}}|\,\D|\mu_m|+	\frac{1}{2}\int_0^{t\wedge\mathfrak{t}_R} \int_{\T^N}\D|{\mu_e}|\leq c(R)	\int_0^{t\wedge\mathfrak{t}_R} \mathcal{D}(r)\,\D r.
	\end{aligned}
\end{equation}
Now we estimate the integral \begin{align*}J_1&\coloneqq\int_{\T^N}\left\lvert\frac{\Div\mathbb{S}\left(\nabla_x\bar{{\bfu}}\right)}{\bar{\varrho}}\right\rvert\langle\nu_{(t,x,\omega)};|\varrho-\bar{\varrho}||\bar{\vu}-\vu|\rangle\,\dx\\
	&\leq c(R)\int_{\T^N}\langle\nu_{(t,x,\omega)};|\varrho-\bar{\varrho}||\bar{\vu}-\vu|\rangle\,\dx.
\end{align*}
Let us define here the cut-off function\begin{equation}\begin{split}\phi_R^\delta&\in C_c^\infty(0,\infty),\,0\leq\phi_R^\delta\leq1 \text{ s.t. } \phi_R^\delta(r)=1 \text{ for } \frac{\underline{\varrho}_R}2\le r\le2\overline{\varrho}_R,\\ &\hspace{15mm} \text{ and }\supp\phi_R^\delta=\left[\frac{\underline{\varrho}_R}2-\delta,2\overline{\varrho}_R+\delta\right].\end{split}\label{cutoff}\end{equation}
Then using $\frac12\leq a^2+(1-a)^2$, we get \begin{align*}
	J_1&\leq c(R)\intTorN{\left\langle\nu_{(t,x,\omega)};|\varrho-\bar{\varrho}||\vu-\bar{\vu}|[\phi_R^\delta(\varrho)]^2\right\rangle}\\
	&\ \ \ +c(R)\intTorN{\left\langle\nu_{(t,x,\omega)};|\varrho-\bar{\varrho}||\vu-\bar{\vu}|[1-\phi_R^\delta(\varrho)]^2\right\rangle}\\
	&\eqqcolon J_{11}+J_{12}.
\end{align*}
Using $ab\leq \frac{a^2+b^2}{2}$, \eqref{inequality0} and \eqref{cutoff}, we get 
\begin{align*}
	J_{11}&\leq c(R)\intTorN{\left\langle\nu_{(t,x,\omega)};|\varrho-\bar{\varrho}|^2[\phi_R^\delta(\varrho)]^2\right\rangle}+c(R)\intTorN{\left\langle\nu_{(t,x,\omega)};|\vu-\bar{\vu}|^2[\phi_R^\delta(\varrho)]^2\right\rangle}\\
	&\leq c(R)\intTorN{\left\langle\nu_{(t,x,\omega)};H(\varrho,\bar{\varrho})+\left(\chi_{\left\{ \frac{\underline{\varrho}_R}2-\delta\leq\varrho\leq\frac{\underline{\varrho}_R}2\right\}}+\chi_{\left\{ {2\overline{\varrho}_R}\leq\varrho\leq{2\overline{\varrho}_R}+\delta\right\}}\right)|\varrho-\bar{\varrho}|^2\right\rangle}\\
	&\ \ \ +c(R)\intTorN{\left\langle\nu_{(t,x,\omega)};|\vu-\bar{\vu}|^2[\phi_R^\delta(\varrho)]^2\right\rangle}\\
	&\leq c(R)\intTorN{\left\langle\nu_{(t,x,\omega)};H(\varrho,\bar{\varrho})+\left[\left(\overline{\varrho}_R-\frac{\underline{\varrho}_R}2+\delta\right)^2+\left(2\overline{\varrho}_R-{\underline{\varrho}_R}+\delta\right)^2\right]\mathcal{O}(\delta)\right\rangle}\\
	&\ \ \ +c(R)\intTorN{\left\langle\nu_{(t,x,\omega)};\chi_{\{ \varrho \geq\frac{\underline{\varrho}_R}2-\delta\}}|\vu-\bar{\vu}|^2\right\rangle}\\
	&\leq c(R)\intTorN{\left\langle\nu_{(t,x,\omega)};H(\varrho,\bar{\varrho})\right\rangle}+\mathcal{O}(\delta)+c(R)\intTorN{\frac{2}{\frac{\underline{\varrho}_R}2-\delta}\left\langle\nu_{(t,x,\omega)};\frac\varrho2|\vu-\bar{\vu}|^2\right\rangle}\\
	&\leq \left(c(R)+\mathcal{O}(\delta)\right)\intTorN{\left\langle\nu_{(t,x,\omega)};H(\varrho,\bar{\varrho})+\frac\varrho2|\vu-\bar{\vu}|^2\right\rangle}+\mathcal{O}(\delta).
\end{align*}
Therefore, letting $\delta\to0$, we get\begin{equation}\label{ineq1}
	J_{11}\leq c(R)	\,\mathcal{E}_{\mathrm{mv}} \big(\varrho ,{\mathbf{m}}\left\vert \right. \bar{\varrho}, \bar{\mathbf{u}}  \big).
\end{equation}
To estimate the integral $J_{12}$, using \eqref{cutoff}, we get
\begin{align*}
	J_{12}&=c(R)\intTorN{\left\langle\nu_{(t,x,\omega)};|\varrho-\bar{\varrho}||\vu-\bar{\vu}|[1-\phi_R^\delta(\varrho)]^2\right\rangle}\\
	&\leq c(R)\intTorN{\left\langle\nu_{(t,x,\omega)};\left(\chi_{\left\{\varrho \leq\frac{\underline{\varrho}_R}{2}\right\}}+\chi_{\left\{\varrho \geq2\overline{\varrho}_R\right\}}\right)\varrho|\vu-\bar{\vu}|\right\rangle}.
\end{align*}
Using $ab\leq \frac{a+ab^2}{2},\,\varrho<1+\varrho^\gamma$, and \eqref{inequality0} we get
\begin{align}\label{ineq2}
	\begin{aligned}
	J_{12}&\leq {c(R)}\intTorN{\left\langle\nu_{(t,x,\omega)};\chi_{\left\{\varrho\leq\frac{\underline{\varrho}_R}{2}\right\}}\frac{\varrho+\varrho|\vu-\bar{\vu}|^2}{2}\right\rangle}\\&\hspace{20mm}+c(R)\intTorN{\left\langle\nu_{(t,x,\omega)};\chi_{\left\{\varrho\geq2\overline{\varrho}_R\right\}}\frac{\varrho+\varrho|\vu-\bar{\vu}|^2}{2}\right\rangle}\\
	&\leq{c(R)}\intTorN{\left\langle\nu_{(t,x,\omega)};\chi_{\left\{\varrho\leq\frac{\underline{\varrho}_R}{2}\right\}}\frac{1+\varrho^\gamma+\varrho|\vu-\bar{\vu}|^2}{2}\right\rangle}\\&\hspace{20mm}+c(R)\intTorN{\left\langle\nu_{(t,x,\omega)};\chi_{\left\{\varrho \geq2\overline{\varrho}_R\right\}}\frac{1+\varrho^\gamma+\varrho|\vu-\bar{\vu}|^2}{2}\right\rangle}\\
	&\leq c(R)\,\mathcal{E}_{\mathrm{mv}} \big(\varrho ,{\mathbf{m}}\left\vert \right. \bar{\varrho}, \bar{\mathbf{u}}  \big).
\end{aligned}
\end{align}
and consequently \eqref{ineq1}, \eqref{ineq2} together give \begin{equation}\label{inequality1}
	J_{1}\leq c(R)	\,\mathcal{E}_{\mathrm{mv}} \big(\varrho ,{\mathbf{m}}\left\vert \right. \bar{\varrho}, \bar{\mathbf{u}}  \big).
\end{equation}

Collecting all the estimates {from} (\ref{inequality1}, \ref{inequality2}, \ref{inequality3}, \ref{inequality6}), we have shown that
\begin{equation}
\begin{aligned}
\label{r4est}
\left\vert\int_0^{t \wedge \mathfrak{t}_R}\mathcal{R}_{\mathrm{mv}} \left(\varrho,{\bf m} \,\big\vert\,\bar{\varrho}, \bar{\mathbf{u}}  \right)(s) \,\D s\right\rvert
\leq c(R)\,\int_0^{t \wedge \mathfrak{t}_R}\mathcal{E}_{\mathrm{mv}} \left(\varrho,{\bf m} \,\big\vert\, \bar{\varrho}, \bar{\mathbf{u}}  \right) 
(s) \,\mathrm{d}s .
\end{aligned}
\end{equation}
Notice that the integrand $\left(\mathbb{S}\left(\nabla_x{\langle\nu_{(s,x,\omega)},\bf u\rangle}\right)-\mathbb{S}(\nabla_x \bar{{\bf u}})\right):\left(\nabla_x{\langle\nu_{(s,x,\omega)},\bf u\rangle}-\nabla_x\bar{\bf u}\right)\geq 0$ for a.e. $(s,x)\in[0,T]\times\T^N$, $\mathbb{P}$-a.s., thus discarding this integral in \eqref{relativeEntropy1} and combining with \eqref{r4est}, and applying Grönwall's inequality on the function $\mathbb{E}[\mathcal{E}_{\mathrm{mv}}]$ we arrive at
\begin{align}
\nonumber
\mathbb{E}\,  \Big[\mathcal{E}_{\mathrm{mv}} \big(\varrho ,{\mathbf{m}}\left\vert \right. \bar{\varrho}, \bar{\mathbf{u}}  \big)  
(t \wedge \mathfrak{t}_R)\Big]\leq c(R,T)\,
\mathbb{E}\,\Big[\mathcal{E}_{\mathrm{mv}} \big(\varrho ,{\mathbf{m}}\left\vert \right. \bar{\varrho}, \bar{\mathbf{u}}  \big)(0)\Big].
\label{relativeEntropy2}
\end{align}
Note that we have
\begin{equation*}
\begin{aligned}
\mathcal{E}_{\mathrm{mv}} \big(\varrho ,{\mathbf{m}}\left\vert \right. \bar{\varrho}, \bar{\mathbf{u}}  \big)  (0) 
&=
\int_{\T^N} \Big \langle \nu_{(0,x,\omega)} ; \frac{1}{2}\varrho_{0}\big\vert \mathbf{u}_{0} - \bar{\mathbf{u}}_{0} \big\vert^2  + H\big(\varrho_{0} ,\bar{\varrho}_0 \big) \Big \rangle \,\mathrm{d}x
\end{aligned}
\end{equation*}
which is zero in expectation by assumption and since $\mathcal{E}_{\mathrm{mv}},\mathcal{D} \geq 0$, this implies $\mathbb{P}$-a.s. $\mathcal{D}(\tau\wedge \mathfrak{t}_R)=0,\,\mathcal{E}_{\mathrm{mv}}(\tau\wedge \mathfrak{t}_R)=0$ for a.e. $\tau \in [0,T]$ and therefore the claim follows.\end{proof}

\section{Incompressible--Inviscid Limit}\label{section incompressible--inviscid limit}
We discuss another application of relative energy inequality in this section to rigorously justify the incompressible --inviscid limit of the system \eqref{stochastic ns density equation}--\eqref{stress tensor}. For this purpose, we first rescale the deterministic counterpart of the stochastic compressible Navier--Stokes system by non-dimensionalization and after combining the terms appropriately and adding a stochastic force term, we get the following system
\begin{align}
	&\D\varrho_\varepsilon+\divv_x\left(\varrho_\varepsilon\vu_\varepsilon\right)=0,\label{M1}\\&\D\left(\varrho_\varepsilon\vu_\varepsilon\right)+\left[\divv_x\left(\varrho_\varepsilon\vu_\varepsilon\otimes\vu_\varepsilon\right)+\frac{1}{\varepsilon^2}\nabla_x p\left(\varrho_\varepsilon\right)\right]\dt=\divv_x\mathbb{S}^\varepsilon\left(\nabla_x\vu_\varepsilon\right)\,\dt+\mathbb{G}\left(\varrho_\varepsilon,\varrho_\varepsilon\vu_\varepsilon\right)\D W,\label{M2}\\& \mathbb{S}^\varepsilon(\Grad \vu_\varepsilon) = \nu_\varepsilon \left( \Grad \vu_\varepsilon + \Grad^t \vu_\varepsilon - \frac{2}{N} \Div \vu_\varepsilon \mathbb{I} \right) + \lambda_\varepsilon \,\Div \vu_\varepsilon \mathbb{I}, \qquad \nu_\varepsilon > 0, \ \lambda_\varepsilon \geq 0,\label{M3}
\end{align} where \[\nu_\varepsilon,\,\lambda_\varepsilon\to0\hspace{2mm}\text{as}\hspace{2mm}\varepsilon\to0.\hspace{2mm}\text{Define}\hspace{2mm}\eta_\varepsilon=\lambda_\varepsilon+\frac{(N-2)\nu_\varepsilon}{N},\hspace{2mm}{\forall N\ge2.}\]

Here we assume a very simple form of $\mathbb{G}$, namely that it is an affine function of density and momentum.\begin{equation}
	\mathbb{G}\left(\varrho,\varrho\vu\right)=\varrho\mathbb{K}+\varrho\vu\mathbb{L}\hspace{3mm}\text{where}\hspace{3mm}\mathbb{K}=\left(K_i\right)_{i\geq1},\, \mathbb{L}=\left(L_i\right)_{i\geq1}\text{ such that }\sum_{i \geq 1}\lvert K_i\rvert+\lvert{L_i}\rvert<\infty.\label{afine}
\end{equation}

Here $K_i,L_i$ are real numbers satisfying the conditions in \eqref{afine} and $\mathbb{K},\mathbb{L}$ are suitable Hilbert-Schmidt operators, so we have a special case of assumptions \eqref{noise condition 1}--\eqref{noise condition 2}. The main advantage of choosing such a $\mathbb{G}$ can be seen below in \eqref{pre-proj}.

The scaling parameter $\varepsilon$ in \eqref{M2} is proportional to the Mach number (the ratio of the characteristic flow velocity and the speed of sound). When we look at a large time scale, from a physical point of view the fluid should behave like an incompressible one if the density is close to $1$ and the velocity is small. For low Mach number (i.e., for small $\varepsilon$) in \eqref{M2}, see \cite{KBNMCS}. Under these circumstances, when we consider the asymptotic limit of solutions $(\varrho_\varepsilon,\vu_\varepsilon)$ for $\varepsilon\to0$, the motion of the fluid is governed by the stochastic incompressible Euler system \begin{align}
	&\D{\bf v}+\left[\Div\left({\bf v}\otimes{\bf v}\right)+\nabla_x\Pi\right]\dt=\Psi({\bf v})\,\D W\label{icM1},\\&\Div{\bf v}=0.\label{icM2}
\end{align}

Here, $\Pi$ is the associated pressure and $\Psi({\bf v})=\mathcal{P}_H\mathbb{G}(1,{\bf v})=\mathbb{G}(1,{\bf v})$, with $\mathcal{P}_H=\operatorname{Id}-\nabla_x\Delta^{-1}\Div$ being the Helmholtz projection onto the space of ``solenoidal'' (divergence free) vector fields.

Our approach is based on a comparison between dissipative measure valued martingale solutions of the primitive system and a local strong solution of the limit system via an application of the relative entropy. Let $\left[ \left(\Omega,\mathfrak{F}, (\mathfrak{F}_{\varepsilon,t})_{t\geq0},\mathbb{P} \right); \nu^\varepsilon_{(t,x,\omega)}, W \right]$ be dissipative measure valued martingale solutions to the system \eqref{M1}--\eqref{M2}, in the sense of Definition \ref{dissipative measure valued martingale solution}. To rigorously justify the process of passing the limit $\varepsilon\to 0$ in \eqref{M1}--\eqref{M2}, we need to make use of the relative energy inequality \eqref{relativeEntropy}. However, since the filtration corresponding to dissipative solutions depend explicitly on $\varepsilon$, we need to justify the usage of the relative energy inequality to the pair $\nu^\varepsilon_{(t,x,\omega)}$ and $\bfv$, for all $\varepsilon>0$. Indeed, this is possible as the strong solution $\bfv$ can be constructed on any given stochastic basis and is adapted to the Brownian motion which is independent of $\varepsilon$. More specifically, for every filtration $\left(\mathfrak{F}_{\varepsilon,t}\right)_{t\geq0}$, the strong solution to the incompressible {Euler} system \eqref{icM1}--\eqref{icM2} is adapted to that filtration, as it is adapted to the Brownian motion and the Brownian motion is adapted to $\left(\mathfrak{F}_{\varepsilon,t}\right)_{t\geq0}$.

\subsection{Strong Solution of the Incompressible Euler System}

Lets assume that we are given the stochastic basis $\left(\Omega,\mathfrak{F},\left(\mathfrak{F}_t\right)_{t\geq0},\mathbb{P}\right)$ and the Wiener process $W$ as identified in the beginning of this paper.

\begin{Definition}\label{def-ssNS}
	Let $\left(\Omega,\mathfrak{F},\left(\mathfrak{F}_t\right)_{t\geq0},\mathbb{P}\right)$ be a stochastic basis with a complete right-continuous filtration, let $W$ be an $\left(\mathfrak{F}_t\right)_{t\geq0}$-cylindrical Wiener process and $\bfv_0$ is a divergence free $\mathfrak{F}_0$-measurable stochastic process. We say $\left(\bfv,\mathfrak{t}\right)$ is a local strong path-wise solution to the incompressible Euler system \eqref{icM1}--\eqref{icM2}, provided \begin{itemize}
		\item $\mathfrak{t}$ is a strictly positive stopping time, for $m>1+\frac{N}2$; $\bfv\in C\left([0,T];W^{m,2}\left(\T^N;\R^N\right)\right)$ $\mathbb{P}$-a.s. is an $\{\mathfrak{F}_t\}$-adapted, predictable stochastic process and \[\expe{\sup_{0\le t \le T}\|\bfv(t,\cdot)\|^p_{W^{m,2}(\T^N;\R^N)}}<\infty\hspace{2mm}\text{for all}\hspace{2mm}1\leq p<\infty.\]
		\item $\mathbb{P}$-a.s., there holds\begin{align}
			&\bfv\left(t\wedge\mathfrak{t}\right)=\bfv_0-\int_0^{t\wedge\mathfrak{t}}\mathcal{P}_H\left[(\bfv\cdot\nabla_x)\bfv\right]\,\D s+\int_0^{t\wedge\mathfrak{t}}\mathbb{G}(1,\bfv)\,\D W_s,\label{ssic}\\
			&\Div \bfv=0,
		\end{align}
	a.e. $(t,x)\in(0,T)\times\T^N$.
	\end{itemize}
\end{Definition}

Regarding the local-in-time existence of strong solutions to the stochastic incompressible {Euler} system \eqref{icM1}--\eqref{icM2}, under certain restrictions imposed on the forcing coefficient $\mathbb{G}$, we refer to the work by Glatt-Holtz and Vicol \cite[Theorem 4.3]{GHVic}. Here, as mentioned before, we assumed a very simple form of $\mathbb{G}$, namely as in \eqref{afine}. One of the crucial benefits of such a choice is that the pressure $\Pi$ can be computed explicitly from \eqref{ssic} and does not contain a stochastic component (in full generality an additional stochastic integral is part of the pressure, see \cite[Section 2]{Breit2015ExistenceTF}). Indeed we note that \begin{equation}\label{pre-proj}
	\nabla_x\Pi=\mathcal{P}_H\left[(\bfv\cdot\nabla_x)\bfv\right]-(\bfv\cdot\nabla_x)\bfv=-\nabla_x\Delta^{-1}\Div\left[(\bfv\cdot\nabla_x)\bfv\right].
\end{equation}
Accordingly the equation \eqref{ssic} reads \begin{equation}\label{ssic1}
	\bfv\left(t\wedge\mathfrak{t}\right)=\bfv_0+\int_0^{t\wedge\mathfrak{t}}\left[-\nabla_x\Pi-(\bfv\cdot\nabla_x)\bfv\right]\,\D s+\int_0^{t\wedge\mathfrak{t}}\mathbb{G}(1,\bfv)\,\D W_s.
\end{equation}

\subsection{Main Result}
We now state the main result of this section related to the rescaled stochastic compressible Navier--Stokes system \eqref{M1}--\eqref{M3}. 

\begin{Theorem}\label{singularlimit}
Let $\tn{G}$ be given as in \eqref{afine}, and the initial data $\vr_{0,\varepsilon}$, $(\vr \vu)_{0,\varepsilon}$, and $\vc{v}_0$ be given such that $\p$-a.s.
	\[
	\left \lbrace \vr_{0,\varepsilon}, (\vr \vu)_{0,\varepsilon}  \in L^\gamma(\tor) \times L^{\frac{2 \gamma}{\gamma + 1}}(\tor;
	\mathbb{R}^N) \ \Bigg| \ \substack{\vr_{0,\varepsilon} \geq \underline \vr > 0,\\ \frac{|\vr_{0,\varepsilon} - 1|}{\varepsilon} \leq \delta(\varepsilon),\\|(\vr \vu)_{0,\varepsilon} - \vc{v}_0 |\leq \delta (\varepsilon)} \right \rbrace,
	\]
	where
	\[
	\delta (\varepsilon) \to 0 \ \mbox{as}\ \varepsilon \to 0,
	\]
	and $\vc{v}_0$ is an $\mathfrak{F}_0$-measurable random variable {satisfying},
	\begin{align*}
		&\vc{v}_0 \in W^{m,2}(\tor;\R^N), \ \Div \vc{v}_0 = 0, \quad \text{$\mathbb{P}$-{\normalfont a.s.}},\\
		&\E{ \| \vc{v}_0 \|_{W^{m,2}(\tor; \R^N)}^p } < \infty, \,\, \mbox{for all}\,\, 1 \leq p < \infty,\,m>1+\frac{N}{2}.
	\end{align*}

	Let $\big[ \big(\Omega,\mathfrak{F}, (\mathfrak{F}_{\varepsilon,t})_{t\geq0},\mathbb{P} \big); \nu^\varepsilon_{(t,x,\omega)}, W \big]$ be a family of dissipative measure valued martingale solutions to the system \eqref{M1}--\eqref{M3} with \[\nu_\varepsilon>0,\,\lambda_\varepsilon\geq0,\,\nu_\varepsilon\to0,\,\lambda_\varepsilon\to0\hspace{2mm}\text{as}\hspace{2mm}\varepsilon\to0\]satisfying the compatibility condition \eqref{bound on defect measures}, and suppose $(\vc{v}, \mathfrak{t})$ defined on the same probability space $\big(\Omega,\mathfrak{F},(\mathfrak{F}_{\varepsilon,t})_{t\geq0},\mathbb{P} \big)$ is a unique local strong solution of the Euler system 
	\eqref{icM1}--\eqref{icM2} driven by the same cylindrical Wiener process $W$ in the sense of Definition~{\normalfont\ref{def-ssNS}}. Then as $\varepsilon \to 0$, $\p$-{\normalfont a.s.}
	\[
	\mathcal{D}^\varepsilon \to 0 \ \mbox{in}\ L^\infty(0,T),
	\]
	\[
	\esssup_{t \in (0,T)} \E \intTorN{
		\left\langle {\nu^\varepsilon_{(t,x,\omega)}}; \frac{1}{2} \vr  \left| \frac{\vc{m}}{\vr} - \vc{v}(t,x) \right|^2 + \frac{P(\vr) - P'(1)(\vr - 1) - P(1)}{\varepsilon^2}  \right\rangle}\,(t \wedge \tau)  \to 0.
	\]
\end{Theorem}
\vspace{2mm}
\begin{Remark}
	We remark that the above theorem asserts that the probability measures ${\nu^\varepsilon_{(t,x,\omega)}}$ shrink to their expected value as $\varepsilon \to 0$, where the latter are characterized
	by the constant value $1$ for the density and the solution $\vc{v}$ of the incompressible system. The situation considered in the above theorem corresponds to the so-called well-prepared data. However, one can extend these results to the case of ill-prepared data (refer to Masmoudi \cite{Masmoudi}, for the related deterministic results) with additional technical computations.
	
\end{Remark}

\subsection[Proof of Theorem 5.2]{Proof of Theorem \ref{singularlimit}}
Let us assume first that $\vv$, with a stopping time $\mathfrak{t}$, is a local strong solution of the stochastic Euler system \eqref{icM1}--\eqref{icM2}. For each $M>0$, let us define \[\tau_M\coloneqq\inf\bigg\{t\in[0,T]:\|\nabla_x\vv(t)\|_{L^\infty(\T^N;\R^N)}>M\bigg\}\] be another stopping time. Thanks to the existence theorem \cite[Theorem 4.3]{GHVic}, we assume, without loss of generality, that $\tau_M\leq\mathfrak{t}$. Following \eqref{refnse}, for $\nu^\varepsilon_{(t,x,\omega)}$-- a dissipative measure valued martingale solution of the rescaled system \eqref{M1}--\eqref{M2}, we define for a.e. $t\in[0,T]$,
\[\begin{split}
\mathcal{E}_{\mathrm{mv}}^{\varepsilon,1} \left( \vr, \vc{m} \,\Big|\, 1, \vc{v} \right)(t)
\coloneqq&\int\limits_{\T^N}
	\left\langle {\nu^{\varepsilon}_{(t,x,\omega)}}; \frac{\vr}{2}   \left| \frac{\vc{m}}{\vr} - \vc{v}(t,x) \right|^2 + \frac{P(\vr) - P'(1)(\vr - 1) - P(1)}{\varepsilon^2}  \right\rangle\dx \\&\quad+ \mathcal{D}^{\varepsilon}(t)\end{split}
\] {to be} the relative energy functional associated to the pair $(1,\vv)$.

Similarly we define  $\mathcal{E}_{\mathrm{mv}}^{\varepsilon,2} \left( \vr, \vc{m} \,\Big|\, 1, \vc{v} \right)$, and $\mathcal{E}_{\mathrm{mv}}^{\varepsilon} \left( \vr, \vc{m} \,\Big|\, 1, \vc{v} \right)$ appropriately.

\subsubsection{Relative Energy Inequality}
As the pair $(r(t),\vU(t))=(1,\vv)(t\wedge\tau_M)$ has the required regularity, it can be used as a test function in the relative energy inequality \eqref{relativeEntropy}. Therefore, by noting that the terms involving $D^d_tr,\,\mathbb{D}^s_tr,\nabla_xP'(r)$ vanish as $r=1$ is a constant, we conclude using $\divv_x\bfv=0,\,D^d_t\bfv=-\nabla_x\Pi-\bfv\cdot\nabla_x\bfv,\,\mathbb{D}^s_t\bfv=\mathbb{G}(1,\bfv)$, that \begin{equation}\label{relEnineq}
	\begin{aligned}
	&\mathcal{E}_{\mathrm{mv}}^{\varepsilon} \left( \vr, \vc{m} \,\Big|\, 1, \vc{v} \right)(t\wedge\tau_M)\\
	&\hspace{5mm}+\int_0^{t\wedge\tau_M}\int_{\T^N}\left(\mathbb{S}^\varepsilon\left(\nabla_x{\langle\nu^\varepsilon_{(s,x,\omega)},\bf u\rangle}\right)-\mathbb{S}^\varepsilon(\nabla_x\bfv)\right):\left(\nabla_x{\langle\nu^\varepsilon_{(s,x,\omega)},\bf u\rangle}-\nabla_x\bfv\right)\dx\,\D s  \\
	&\leq\intTorN{
		\left\langle {\nu^{\varepsilon}_{(0,x,\omega)}}; \frac{1}{2} \vr  \left| \frac{\vc{m}}{\vr} - \vc{v}(0,x) \right|^2 + \frac{1}{\varepsilon^2} \Big( P(\vr) - P'(1)(\vr - 1) - P(1) \Big) \right\rangle}\\
	&\hspace{5mm}+M_{RE}(t\wedge\tau_M)\\
	&\hspace{5mm}+\int_0^{t\wedge\tau_M}\int_{\T^N}\mathbb{S}^\varepsilon\left(\nabla_x{\bfv}\right):\left(\nabla_x{\bfv}-\nabla_x\left\langle\nu^\varepsilon_{{(s,x,\omega)}},\bfu\right\rangle\right)\dx\,\D s
	\\&\hspace{5mm}+\int_0^{t\wedge\tau_M}\int_{\T^N}\langle \nu^\varepsilon_{{(s,x,\omega)}}; \varrho \mathbf{v} - {\bf m} \rangle \cdot[-\nabla_x\Pi] \,\dx\,\D s \\
	&\hspace{5mm}+\int_0^{t\wedge\tau_M}\int_{\T^N} \left\langle \nu^\varepsilon_{(s,x,\omega)}; \frac{({\bf m}-\varrho  \mathbf{v})\otimes (\varrho  \mathbf{v}-{\bf m}) }{\varrho} \right\rangle: \nabla_x {\bfv} \,\dx\,\D s
	\\
	&\hspace{5mm}+
	\frac{1}{2}
	\sum_{k\in\mathbb{N}}\int_0^{t\wedge\tau_M}
	\int_{\T^N} \left \langle \nu^\varepsilon_{(s,x,\omega)}; \varrho\left\vert \frac{\mathbf{G}_k(\varrho,{\bf m})}{\varrho}  -\mathbf{G}_k(1,\bfv)  \right\vert^2\right\rangle \,\dx\,\D s \\
	&\hspace{5mm}-\int_0^{t\wedge\tau_M}  \int_{\T^N} \nabla_x \mathbf{v} : \D\mu_m + \frac12\int_0^{t\wedge\tau_M}\int_{\T^N} \D \mu_e.
	\end{aligned}
\end{equation}

Similarly to the proof of Theorem \ref{weakstrong uniqueness for dmv solution nse}, we show that the terms on the right-hand side of \eqref{relEnineq} can be absorbed by means of a Grönwall-type argument. To see this, we first observe that,
\begin{equation}\label{estimate1}
	\begin{aligned}
		&\int_{\T^N}\mathbb{S}^\varepsilon\left(\nabla_x{\bfv}\right):\left(\nabla_x{\bfv}-\nabla_x\left\langle\nu^\varepsilon_{{(s,x,\omega)}},\bfu\right\rangle\right)\dx\\
		&=\frac12\int_{\T^N}\left(\mathbb{S}^\varepsilon\left(\nabla_x{\langle\nu^\varepsilon_{(s,x,\omega)},\bf u\rangle}\right)-\mathbb{S}^\varepsilon(\nabla_x\bfv)\right):\left(\nabla_x{\langle\nu^\varepsilon_{(s,x,\omega)},\bf u\rangle}-\nabla_x\bfv\right)\dx\\
		&\hspace{5mm}+\frac12\int_{\T^N}\left(\mathbb{S}^\varepsilon\left(\nabla_x{\bfv}\right)+\mathbb{S}^\varepsilon\left(\nabla_x{\langle\nu^\varepsilon_{(s,x,\omega)},\bf u\rangle}\right)\right):\left(\nabla_x{\bfv}-\nabla_x\left\langle\nu^\varepsilon_{{(s,x,\omega)}},\bfu\right\rangle\right)\dx,
	\end{aligned}
\end{equation} and\footnote{{We remark that the estimate \eqref{estimate2} holds true when $N=1$.}} \begin{equation}\label{estimate2}
\begin{aligned}
	&\left(\mathbb{S}^\varepsilon\left(\nabla_x{\bfv}\right)+\mathbb{S}^\varepsilon\left(\nabla_x{\langle\nu^\varepsilon_{(s,x,\omega)},\bf u\rangle}\right)\right):\left(\nabla_x{\bfv}-\nabla_x\left\langle\nu^\varepsilon_{{(s,x,\omega)}},\bfu\right\rangle\right)\\
	&\quad=\nu_\varepsilon|\nabla_x\bfv|^2+\eta_\varepsilon|\divv_x\bfv|^2-\nu_\varepsilon|\nabla_x\langle\nu^\varepsilon_{(s,x,\omega)};\bfu\rangle|^2-\eta_\varepsilon|\divv_x\langle\nu^\varepsilon_{(s,x,\omega)};\bfu\rangle|^2\\
	&\quad\quad\leq c\,(\nu_\varepsilon+\eta_\varepsilon) M^2\leq c(N)M^2(\nu_\varepsilon+\lambda_\varepsilon).
\end{aligned}
\end{equation}

Next we see, \begin{equation}\label{estimate3}
\begin{aligned}
&\left\vert\int_{\T^N} \left\langle \nu^\varepsilon_{(s,x,\omega)}; \frac{({\bf m}-\varrho  \mathbf{v})\otimes (\varrho  \mathbf{v}-{\bf m}) }{\varrho} \right\rangle: \nabla_x {\bfv} \,\dx\right\vert\\
&\qquad\leq c(M)\int_{\T^N} \left\langle \nu^\varepsilon_{(s,x,\omega)}; \frac{\left\vert({\bf m}-\varrho  \mathbf{v})\right\vert^2}{2\varrho} \right\rangle \,\dx\\
&\qquad\leq c(M)\,\mathcal{E}_{\mathrm{mv}}^{\varepsilon} \left( \vr, \vc{m} \,\Big|\, 1, \vc{v} \right)(s).
\end{aligned}
\end{equation}
Similarly, using \eqref{afine}, we get\begin{equation}\label{estimate4}
	\begin{aligned}
		&\sum_{k\in\mathbb{N}}\int_{\T^N} \left \langle \nu^\varepsilon_{(s,x,\omega)}; \varrho\left\vert \frac{\mathbf{G}_k(\varrho,{\bf m})}{\varrho}  -\mathbf{G}_k(1,\bfv)  \right\vert^2\right\rangle \,\dx\\&\qquad\leq \sum_{k \in\mathbb{N}}\int_{\T^N}\left\langle\nu^\varepsilon_{(s,x,\omega)};\frac{1}{\varrho}\left\vert{\bf m}-\varrho\bfv\right\vert^2L_k^2\right\rangle\,\dx\leq c\,\mathcal{E}_{\mathrm{mv}}^{\varepsilon} \left( \vr, \vc{m} \,\Big|\, 1, \vc{v} \right)(s),
	\end{aligned}
\end{equation}
and in addition, since the compatibility condition \eqref{bound on defect measures} is satisfied for all $\varepsilon$, we deduce
\begin{equation} \label{estimate5}
	\frac12 \int_0^{t\wedge \tau_M} \int_{\T^N} \D \mu^{\varepsilon}_e - \int_0^{t \wedge \tau_M} \int_{\mathbb{T}^N} \Grad \vc{v}:{\rm d}\mu^{\varepsilon}_m 
	\leq c\int_0^{t \wedge \tau_M} \mathcal{D}^\varepsilon(s)\,\D s.
\end{equation}
As the initial data is well-prepared, we get
\begin{equation} \label{estimate6}
	\begin{split}
	&\E \left[\intTorN{
		\left\langle\nu^{\varepsilon}_{(0,x,\omega)}; \frac{1}{2} \vr  \left\vert \frac{\vc{m}}{\vr} - \vc{v}_0(x) \right\vert^2 + \frac{P(\vr) - P'(1)(\vr - 1) - P(1)}{\varepsilon^2} \right\rangle} \right] \le \delta(\varepsilon),\\&\hspace{4cm}\text{where}\hspace{3mm}\delta(\varepsilon) \to 0 \ \mbox{as}\ \varepsilon \to 0.
	\end{split}
\end{equation} Therefore taking expectation of \eqref{relEnineq} and invoking all the estimates (\ref{estimate1}, \ref{estimate2}, \ref{estimate3}, \ref{estimate4}, \ref{estimate5} and \ref{estimate6}) into it, we get\begin{equation}\label{estimate7}
\begin{aligned}
	&\E\left[\mathcal{E}_{\mathrm{mv}}^{\varepsilon} \left( \vr, \vc{m} \,\Big|\, 1, \vc{v} \right)(t\wedge\tau_M)\right]\\
	&\hspace{5mm}+\frac12\,\E\left[\int_0^{t\wedge\tau_M}\int_{\T^N}\left(\mathbb{S}^\varepsilon\left(\nabla_x{\langle\nu^\varepsilon_{(s,x,\omega)},\bf u\rangle}\right)-\mathbb{S}^\varepsilon(\nabla_x\bfv)\right):\left(\nabla_x{\langle\nu^\varepsilon_{(s,x,\omega)},\bf u\rangle}-\nabla_x\bfv\right)\dx\,\D s\right]  \\
	&\leq\delta(\varepsilon)+ c(N)M^2(\nu_\varepsilon+\lambda_\varepsilon)+c(M)\int_0^{t\wedge\tau_M}\mathcal{E}_{\mathrm{mv}}^{\varepsilon} \left( \vr, \vc{m} \,\Big|\, 1, \vc{v} \right)(s)\,\D s\\
	&\hspace{5mm}+\E\left[\int_0^{t\wedge\tau_M}\int_{\T^N}\langle \nu^\varepsilon_{{(s,x,\omega)}}; \varrho \mathbf{v} - {\bf m} \rangle \cdot[-\nabla_x\Pi] \,\dx\,\D s\right].
\end{aligned}
\end{equation} Now, we write the last integral in \eqref{estimate7} as follows :

By using the incompressibility condition $\divv_x\bfv=0$ and an integration by parts give
\begin{equation}\label{estimate8}
	\begin{aligned}
	\E\int_0^{t\wedge\tau_M}\int_{\T^N}\langle \nu^\varepsilon_{{(s,x,\omega)}}; \varrho \mathbf{v}\rangle \cdot\nabla_x\Pi \,\dx\,\D s=	\varepsilon\,\E\int_0^{t\wedge\tau_M}\int_{\T^N}\left\langle \nu^\varepsilon_{{(s,x,\omega)}}; \frac{\varrho-1}\varepsilon\right\rangle\bfv \cdot\nabla_x\Pi \,\dx\,\D s,
	\end{aligned}
\end{equation} and the density equation \eqref{measure valued density equation} with $\mu_c=0$, gives \begin{equation}\label{estimate9}
\begin{aligned}
			&\E\left[\int_0^{t \wedge \tau_M} \intTorN{ \left< \nu^\varepsilon_{(s,x,\omega)}; \vc{m} \right> \cdot \Grad \Pi }\,\D s\right]  =\E \left[ \intTorN{ \left<\nu^\varepsilon_{(s,x,\omega)}; \vr \right> \Pi } \right]_{s = 0}^{s= t  \wedge \tau_M}\\&\hspace{3cm}=
			\varepsilon\,\E \left[ \intTorN{ \left<\nu^\varepsilon_{(s,x,\omega)}; \frac{\vr - 1}{\varepsilon} \right> \Pi } \right]_{s = 0}^{s = t  \wedge \tau_M}.
\end{aligned}
\end{equation} Therefore, after ignoring the positive integral in the left hand side of \eqref{estimate7}, it becomes \begin{equation}\label{estimate10}
\begin{aligned}
	\E\left[\mathcal{E}_{\mathrm{mv}}^{\varepsilon} \left( \vr, \vc{m} \,\Big|\, 1, \vc{v} \right)(t\wedge\tau_M)\right]&\leq\delta(\varepsilon)+ c(N)M^2(\nu_\varepsilon+\lambda_\varepsilon)\\&\hspace{5mm}+c(M)\,\E\int_0^{t\wedge\tau_M}\mathcal{E}_{\mathrm{mv}}^{\varepsilon} \left( \vr, \vc{m} \,\Big|\, 1, \vc{v} \right)(s)\,\D s\\
	&\hspace{5mm}-\varepsilon\,\E\left[\int_0^{t\wedge\tau_M}\int_{\T^N}\left\langle \nu^\varepsilon_{{(s,x,\omega)}}; \frac{\varrho-1}\varepsilon\right\rangle\bfv \cdot\nabla_x\Pi \,\dx\,\D s\right]\\
	&\hspace{6mm}+\varepsilon\,\E \left[ \intTorN{ \left<\nu^\varepsilon_{(s,x,\omega)}; \frac{\vr - 1}{\varepsilon} \right> \Pi } \right]_{s = 0}^{s = t  \wedge \tau_M}.
\end{aligned}
\end{equation}
In order to control the last two terms of \eqref{estimate10}, we invoke $r=1,{\bf U}=0$ in \eqref{relativeEntropy} and obtain\begin{equation}
\begin{split}
	&\expe{ \intTorN{ \left \langle \nu^\varepsilon_{(t,x,\omega)}; \frac{1}{2} \vr |\vu|^2 + \frac{P(\vr) - P'(1)(\vr - 1) - P(1)}{\varepsilon^2} \right \rangle}\,(t \wedge \tau_M) + \mathcal{D}^{\varepsilon}(t \wedge \tau_M) } \\
	&\qquad \leq
	\expe{ \intTorN{ \left \langle \nu^\varepsilon_{(0,x,\omega)};\frac{1}{2} \vr |\vu|^2 +  \frac{P(\vr) - P'(1)(\vr - 1) - P(1)}{\varepsilon^2}\right \rangle}\,(0) }.
\end{split}
\end{equation} Thus, if the right hand side of the above inequality is bounded uniformly for $\varepsilon\to0$, we get the following uniform bounds (by using \eqref{inequality0} and setting $\gamma_*=\min\lbrace \gamma, 2\rbrace$)
\begin{align}\label{estimate11}
\E\int_{\T^N} \left< \nu^\varepsilon_{(t,x,\omega)}; \frac12 \vr |\vu|^2 \right>\dx\, (t\wedge \tau_M) \le c, \quad
\E \int_{\T^N} \left< \nu^\varepsilon_{(t,x,\omega)}; \frac{|\vr -1|^{\gamma_*}}{\varepsilon^2}\right>\dx\,(t \wedge \tau_M)\le c.
\end{align}
Now using \eqref{pre-proj}, the continuity of $\nabla_x\Delta^{-1}\Div$, regularity of $\vc{v}$, Hölder's inequality, and \eqref{estimate11}, we obtain 
\begin{align}\label{estimate12}
	\begin{aligned}
&\left| \expe{ \int_0^{t \wedge \tau_M} \intTorN{ \left< \nu^\varepsilon_{(s,x,\omega)}; \frac{ \varrho - 1}{\varepsilon} \right>\,\Grad \Pi \cdot {\bfv} }\,\D s } \right| \\
& \qquad \le\varepsilon^{\frac{2-\gamma_*}{\gamma_*}} \left(\E\int_0^{t\wedge\tau_M}\int_{\T^N}\left\langle\nu^\varepsilon_{(s,x,\omega)};\frac{|\varrho-1|^{\gamma_*}}{\varepsilon^2}\right\rangle\,\dx\,\D s\right)^\frac{1}{\gamma_*}\left\|\nabla_x\Pi\cdot\bfv\right\|_{\gamma_*'}\\
&\qquad\leq \mathcal{O}\left(\varepsilon^{\frac{2-\gamma_*}{\gamma_*}}\right)\left\|\nabla_x\Pi\right\|_{2\gamma_*'}\left\|\bfv\right\|_{2\gamma_*'}\\
&\qquad=\mathcal{O}\left(\varepsilon^{\frac{2-\gamma_*}{\gamma_*}}\right)\left\|\bfv\cdot\nabla_x\bfv\right\|_{2\gamma_*'}\left\|\bfv\right\|_{2\gamma_*'}\\
&\qquad\leq\mathcal{O}\left(\varepsilon^{\frac{2-\gamma_*}{\gamma_*}}\right).
\end{aligned}
\end{align}In an exactly similar manner, we also get\begin{equation}\label{estimate13}
\E \left[ \intTorN{ \left<\nu^\varepsilon_{(s,x,\omega)}; \frac{\vr - 1}{\varepsilon} \right> \Pi } \right]_{s = 0}^{s = t  \wedge \tau_M}=\mathcal{O}\left(\varepsilon^{\frac{2-\gamma_*}{\gamma_*}}\right).
\end{equation}
Finally, using \eqref{estimate12}, \eqref{estimate13} in \eqref{estimate10}, we deduce
\begin{align}
	\begin{aligned}
		\E\mathcal{E}_{\mathrm{mv}}^{\varepsilon} \left( \vr, \vc{m} \,\Big|\, 1, \vc{v} \right)(t\wedge\tau_M)&\leq\mathcal{O}\left(\varepsilon^{\frac{2}{\gamma_*}}+\delta(\varepsilon)+\nu_\varepsilon+\lambda_\varepsilon\right)\\&\hspace{5mm}+c(M)\int_0^{t\wedge\tau_M}\E\mathcal{E}_{\mathrm{mv}}^{\varepsilon} \left( \vr, \vc{m} \,\Big|\, 1, \vc{v} \right)(s)\,\D s.
	\end{aligned}
\end{align}
Applying, Grönwall's lemma on the above inequality gives\[\E\mathcal{E}_{\mathrm{mv}}^{\varepsilon} \left( \vr, \vc{m} \,\Big|\, 1, \vc{v} \right)(t\wedge\tau_M)\leq\mathcal{O}\left(\varepsilon^{\frac{2}{\gamma_*}}+\delta(\varepsilon)+\nu_\varepsilon+\lambda_\varepsilon\right){\rm e}^{c(M)T}\to0\hspace{2mm}\text{as}\hspace{2mm}\varepsilon\to0,\] and the result stated in Theorem \ref{singularlimit} follows.

\subsection*{Acknowledgments}
The author expresses their sincere gratitude to the reviewers for their thoughtful comments and constructive criticisms that have contributed significantly to improving the clarity and overall quality of this research.

\subsection*{Disclosure Statement} The author(s) declare no competing interests.

\printbibliography
\end{document}